\documentclass[12pt]{article}
\usepackage{amsmath,amsfonts,amssymb,amsthm}
\usepackage[dvips]{graphicx}
\usepackage{graphics}
\usepackage{subfig,sidecap}
\usepackage{alltt,makeidx,enumerate}
\allowdisplaybreaks
\theoremstyle{plain}
\newtheorem{theorem}{Theorem}[section]
\newtheorem{lemma}[theorem]{Lemma}
\newtheorem{proposition}[theorem]{Proposition}

\theoremstyle{definition}

\newtheorem{remark}[theorem]{Remark}

\newcommand{\me}{\mathrm{e}}
\newcommand{\dif}{\mathrm{d}}

\newcommand{\deru}[1]{#1^\prime}

\newcommand{\tsg}[1]{\boldsymbol{#1}}
\newcommand{\norm}[1]{\Vert#1\Vert}
\newcommand{\abs}[1]{|#1|}
\newcommand{\Real}{\mathbb R}
\newcommand{\eps}{\varepsilon}
\newcommand{\To}{\rightarrow}
\newcommand{\supp}[1]{\mathrm{supp}(#1)}
\newcommand{\limite}[2]{\lim_{#1\rightarrow #2}}
\newcommand{\td}[2]{\frac{\dif#1}{\dif#2}}

\newcommand{\wsconv}[2]{#1\stackrel{*}{\rightharpoonup}#2}

\newcommand{\R}{{\mathbb R}}
\newcommand{\ra}{\rightarrow}
\newcommand{\bu}{{\mathbf u}}

\newcommand{\ts}[1]{\mathbf{#1}}
\newcommand{\mc}[1]{\mathcal{#1}}

\newcommand{\set}[1]{\left\{#1\right\}}
\newcommand{\bx}{{\mathbf x}}
\newcommand{\M}{\mathrm{M}}

\newcommand{\bphi}{{\mathbf {\phi}}}

\newcommand{\Om}{{\Omega}}

\DeclareMathOperator{\sgn}{sgn}

\newenvironment{prueba}{\noindent\textit{Proof}:}{  \qed\\\indent}

\begin{document}

\title{
The repulsion property in nonlinear elasticity and a numerical
scheme to circumvent it}

\author{Pablo V. Negr\'on--Marrero\\
Department of Mathematics\\
   University of Puerto Rico\\
   Humacao, Puerto Rico 00791\and
  Jeyabal Sivaloganathan \\
  Department of Mathematical Sciences \\
  University of Bath \\
  Bath \  BA2 7AY, U.K.}

\maketitle

\begin{abstract}
For problems in the calculus of cariations that exhibit the
Lavrentiev phenomenon, it is known that the \textit{repulsion property} holds,
that is, if one approximates the global minimizer in these problems by smooth
functions, then the approximate energies will blow up. Thus standard numerical
schemes, like the finite element method, may fail when applied directly to these
type of problems. In this paper we prove that the repulsion property holds for
variational problems in three dimensional elasticity that exhibit cavitation.
In addition we propose a numerical scheme that circumvents the repulsion 
property, which is an adaptation of the Modica and Mortola functional for phase 
transitions in liquids, in which the phase function is coupled to the mechanical 
part of the stored energy functional, via the determinant of the deformation 
gradient. We show that the corresponding approximations by this method satisfy 
the lower bound $\Gamma$--convergence property in the multi-dimensional 
non--radial case. The convergence to the actual cavitating minimizer is 
established for a spherical body, in the case of radial deformations, and for 
the case of an elastic fluid without assuming radial symmetry.
\end{abstract}

\vspace{.3in}




\section{Introduction}
One-dimensional problems in the Calculus of Variations
that exhibit the Lavrentiev Phenomenon \cite{La1926} have been well studied 
(see \cite{BaMi}). 
 A typical result in such problems, is that
the infimum of a given integral functional
$$
I(u)=\int_a^b L(x, u(x), u'(x))\,\dif x,
$$
on the admissible set of Sobolev functions
$$
\mathcal{A}_p=\{ u\in W^{1,p} ((a,b))\ | \ u(a)=\alpha, \ u(b)=\beta \}, \ p>1,
$$
is strictly greater than its infimum in the corresponding set of
absolutely continuous functions $$\mathcal{A}_1=\{ u\in W^{1,1} ((a,b))\ | \ 
u(a)=\alpha, \ u(b)=\beta \},
$$
i.e., for $p>1$,  $$
\inf _{u\in \mathcal{A}_1}I(u) < \inf _{u\in \mathcal{A}_p}I(u).
$$
Moreover, it has been shown (see \cite[Theorem 5.5]{BaMi}) in a number  of 
cases 
that if the Lavrentiev Phenomenon occurs, then a``repulsion property" holds 
when 
trying to approximate a minimiser by more regular functions: that if $u_0\in 
\mathcal{A}_1 $ is a minimiser  of $I$ on $\mathcal{A}_1 $, and if 
$(u_n)\subset 
\mathcal{A}_p, \ p>1, $ satisfies $u_n\rightarrow u_0$ almost everywhere,
then  $I(u_n)\rightarrow \infty$ as $n\rightarrow \infty$.  We refer to 
\cite{Fe2007} for results on the repulsion property for multi--dimensional 
problems of the calculus of variations that exhibit the Lavrentiev phenomena. 
In particular, they show for a quite general class of functionals, that for any 
minimizer of a problem that exhibits the Lavrentiev phenomena, 
there exists a sequence of ``smooth'' functions converging (strongly) in 
$W^{1,p}$ to the minimizer (for some $p$), but for which the values of the 
functional on the sequence, tend to infinity.

The Lavrentiev phenomena shows up as well for problems in 
elasticity (cf. \cite{Ba82}, \cite{FoHrMi2003}). In the first part of this 
paper we prove a corresponding repulsion property for variational
problems in elasticity that exhibit cavitation (Theorem \ref{thm:repprop}). 
Our result, although more specific in terms of the functionals than those 
considered in \cite{Fe2007}, but which are relevant to elasticity and for which 
minimizers are guaranteed to exist, is stronger in the sense that the 
function we are trying to approximate need not be a minimizer, just a 
cavitating deformation with finite energy, and the convergence in the repulsion 
property is weakly in $W^{1,p}$. The result, though straightforward  to prove, 
does not appear to have been noted previously and has implications for the 
design of numerical methods to detect cavitation instabilities in nonlinear 
elasticity. In particular, from the proof of Theorem \ref{thm:repprop}, it 
becomes evident that the critical term in the stored energy function, in 
relation to the repulsion property, is the compressibility term (the $h(\cdot)$ 
term in \eqref{SEF}).

The numerical aspects of computing cavitated solutions are challenging due to 
the singular nature of such deformations. The work of Negr\'on--Marrero 
\cite{Ne90} generalized to the multidimensional case of elasticity a 
method introduced by Ball and Knowles \cite{BaKn87} for one dimensional 
problems, which is based on a decoupling technique that detects singular 
minimizers and avoids the Lavrentiev phenomena. The convergence 
result in \cite{Ne90} involved a very strong condition on the adjoints of the 
finite element approximations which among other things excluded cavitated 
solutions. The element removal method introduced by Li (\cite{Li92},
\cite{Li95}) improves upon this by penalizing or excluding the elements of the 
finite element grid where the deformation gradient becomes very large. We refer 
also to the works of Henao and Xu \cite{HX} and Lian and Li (\cite{LiLi2011a}, 
\cite{LiLi2011b}).

Motivated by the result in Theorem \ref{thm:repprop}, we propose in Section 
\ref{sec:num} a numerical scheme for computing cavitating deformations, that 
avoids or works around the repulsion property by using nonsingular or smooth 
approximations. The idea is to introduce a decoupling or phase function on the 
determinant of the competing deformations, together with an extra term in the 
energy functional that forces the phase function to assume either small or very 
large values, and penalizes for the corresponding transition regions. This 
extra term is a variant of the Modica and Mortola functional considered in 
\cite{Mo1987} for phase transitions in liquids, The Modica and Mortola 
functional has the form
\[
 \int_{\Omega}\left[\eps \norm{\nabla
v(\ts{x})}^2+\frac{v^2(\ts{x})(1-v(\ts{x}))^2}{\eps}
\right]\dif\ts{x},
\]
where $\eps>0$ is an approximation parameter. The second term in this 
functional forces the phase function $v$ to take 
values either close to zero or one. On the other hand, the extra term in our 
functional (cf. \eqref{modfunct}) is given, for any $\tau>0$, by
\[
\int_{\Omega}\left[\frac{\eps^\alpha}{\alpha} \norm{\nabla
v(\ts{x})}^\alpha+\frac{1}{q\eps^q}
\phi_\tau(v(\ts{x}))\right]\dif\ts{x},
\]
where $\alpha>1$, $\frac{1}{\alpha}+\frac{1}{q}=1$, and 
$\phi_\tau:\Real\To[0,\infty)$ is a continuous function such that 
$\supp{\phi_\tau}\subset(0,\tau)$. The second term in this functional 
penalizes for regions where the phase function $v$ is positive but less 
than $\tau$, but does not penalize for values of $v$ greater than $\tau$. This 
phase function, which in addition is required to satisfy the constraint 
$v\ge0$, is now coupled to the mechanical energy through the compressibility 
term (cf. \eqref{modfunct}). Our approach is similar in spirit to that of 
Henao, Mora--Corral, and Xu \cite{HeMoXu2015}, but with some major differences.  
In particular, they employ two phase functions $v$ and $w$, with the $v$ similar 
to ours but coupled to the mechanical energy as a factor multiplying the 
original stored energy function, and $w$ defined on the deformed configuration. 
The extra terms are of the Ambrosio--Tortorelli \cite{AmTo1990} type for $v$ and 
of the Modica--Mortola type for $w$. As the parameter $\eps$ goes to zero, these 
extra terms in the energy functional, allow for the approximation of 
deformations that can exhibit cavitation or fracture.

In Theorems \ref{DPconvG} and \ref{reppnumsch} we show 
that the proposed scheme has the lower bound $\Gamma$--convergence property. 
Moreover, the approximating functions, which now depend on both $\tau$ and 
$\eps$, converge weakly (for a subsequence and with $\tau\To\infty$ and 
$\eps\To0^+$) in $W^{1,p}$ with $m-1<p<m$ and $m$ the space dimension, to a 
function $\ts{u}^*$ whose distributional determinant is a positive Radon 
measure. The singular part of this measure characterizes the points or regions 
in the reference configuration, where discontinuities of cavitation or fracture 
type can occur.
Further refinements of these results are discussed for two important cases: for 
radial deformations of an spherical body 
(Section \ref{rad:case}) and for an elastic fluid (Section \ref{elast:case}) 
for a spherical body as well. In the former we show (Theorem \ref{convgRad}) 
that the approximations of the proposed decoupled--penalized method, converge to 
the radial cavitating solution. Also we show that the minimizers of the 
penalized functionals (cf. \eqref{modfunctR}) satisfy the corresponding versions 
of the Euler--Lagrange equations, together with some numerical simulations. For 
the elastic fluid case, using the concept of symmetric decreasing rearrangements 
of functions, we show that the corresponding penalized minimizers must be 
radial. Consequently, using the result for radial deformations of Section 
\ref{rad:case}, we get the convergence of the decoupled--penalized 
approximations to the cavitating minimizer in this case as well. We present 
some numerical simulations without assuming radial symmetry, for a two 
dimensional version of the elastic fluid case.

\section{Background}
Let $\Omega \subset \R ^m$ ($m=2$ or $m=3$) denote the region
occupied by a nonlinearly elastic body in its reference
configuration. A deformation of the body corresponds to a map $\bu
: \Omega \ra \R^m$, $\bu \in W^{1,1}(\Omega )$, that is one-to-one
almost everywhere and satisfies the condition
\begin{equation}
\label{invertibility} \det \nabla \bu (\bx ) > 0 \ \ \mbox{for
a.e. } \bx\in\Omega .
\end{equation}
In hyperelasticity the energy stored under such a deformation is
given by
\begin{equation}\label{energy}
E(\bu ) = \int_{\Omega} \ W (\nabla \bu (\bx ))\,\dif\bx,
\end{equation}
where $W: \M^{m\times m}_+ \ra [0,\infty)$ is
the stored energy function of the material and $\M^{m \times m}
_+$ denotes the set of real $m\times m$ matrices with positive
determinant. We consider the displacement problem in which we
require
\begin{equation}\label{boundary}
\bu (\bx )  = \bu^h (\bx)  \ \mathrm{for}\ \bx \in
\partial \Omega , \quad \bu^h(\bx )\equiv \mathbf{A} \bx \ ,
\end{equation}
where $\mathbf{A} \in \M^{m\times m}_+$ is fixed.
Let $\Omega \subset\subset \Omega^e$, where $\Omega^e $ is a bounded, open, 
connected
set with smooth boundary.

\subsection{The distributional determinant}
\noindent
If $\bu\in W^{1,p}(\Omega )$ satisfies
\eqref{boundary},
then we define its homogeneous extension
$\bu_e:\Omega^e\rightarrow \mathbb{R}^m$ by
\begin{equation}\label{extention}
\bu_e(\bx)= \left\{\begin{array}{ll}
                \bu(\bx) & \mbox{if $\bx\in\Omega$,}\\
               \mathbf{ A}\bx   & \mbox{if $\bx\in\Omega^e\backslash\Omega$,}
                      \end{array}
               \right.
\end{equation}
and note that $\bu_e\in W^{1,p}(\Omega^e ;\mathbb{R}^m)$. For
$p>m^2/(m+1)$,

\begin{equation}\label{DDet}
\mbox{Det}\nabla \bu (\bphi ) := -\int_{\Om} \frac{1}{m}\left(
[\mbox{adj} \nabla \bu] \bu \right)
\cdot \nabla \bphi \, \dif \bx, \ \ \
\forall \bphi \in C^\infty _0 (\Omega),
\end{equation}
\noindent is a well-defined distribution. (Here $\mbox{adj}
\nabla \bu$ denotes the adjugate matrix of $\nabla \bu$, that is,
the transposed matrix of cofactors of $\nabla \bu $.) The
definition follows from the well-known formula for expressing
$\det \nabla \bu$ as a divergence. (See, e.g.,\ \cite{Mu_A} for
further details and references.)

Next suppose that $\bu \in W^{1,p}(\Omega), \ p>m-1$, and that $\bu_e$ satisfies
condition (INV)
(introduced by Muller and Spector in  \cite{MuSp}) on $\Omega^e$. Then $\bu_e
\in L^\infty_{loc}
(\Omega^e )$ and hence $\mathrm{Det}(\nabla\bu)$ is again a well-defined
distribution. Moreover, it follows from \cite[Lemma 8.1]{MuSp}
that if $\bu$ further satisfies $\det\nabla \bu >0$ a.e.\ then
$\mathrm{Det}\nabla \bu$ is a Radon measure and
\begin{eqnarray}
\mathrm{Det}\nabla\bu =(\det\nabla\bu)\mathcal{L}^m+\mu_s,
\end{eqnarray}
\noindent where $\mu_s$ is singular with respect to Lebesgue measure
$\mathcal{L}^m$. We first consider the case when $\mu_s$ is a Dirac
measure\footnote{Other assumptions on the support of the singular
measure $\mu ^s$  may be relevant for modelling different forms of
fracture. See also \cite{Mu_B} for further results on the singular
support of the distributional Jacobian.} of the form $\alpha
\delta _{\bx_0}$ (where $\alpha
>0$ and $\bx_0 \in \Omega$) which corresponds to
$\bu$ creating a cavity\footnote{Note that such a cavity need not be spherical.}
of volume $\alpha$ at the point
$\bx_0$.
Following \cite{SiSp1}, we fix $\bx_0 \in \Omega$ and define the set of 
admissible deformations by
\begin{eqnarray} \label{set}
\mathcal{A}_{\bx_0} =\{\bu\in W^{1,p} (\Omega) :\bu |_{\partial\Omega}=\bu^h,\
\bu_e\ {\rm satisfies}\ {\rm (INV)\ on}\ \Omega, \\ \nonumber
\det\nabla\bu>0\ {\rm a.e.},\ \
\mathrm{Det}\nabla\bu=(\det\nabla\bu)\mathcal{L}^m+\alpha_{\bu}\delta_{\bx_0}
\},
\end{eqnarray}
\noindent where $\alpha_\bu\geq 0$ is a scalar depending
on the map $\bu$, and $\delta_{\bx_0}$
denotes the Dirac measure with support at $\bx_0$.
Thus, $\mathcal{A}_{\bx_0}$ contains maps $\bu$ that produce a cavity of volume
$\alpha_{\bu}$ located at $\bx_0\in\Omega$. We will say that the deformation 
$\bu\in \mathcal{A}_{\bx_0}$ is
\emph{singular} if $\alpha_\bu> 0$.

\section{Singular Minimisers, Deformations and the Repulsion
Property}\label{sec3}
In this note, for simplicity of exposition, we consider stored energy functions 
of the form
\begin{equation}\label{SEF}
W(\mathbf{F})=\tilde W(\mathbf{F}) +h(\det \mathbf{F}) \ \ \mathrm{for \ }
\mathbf{F}\in M^{m\times m}_+,
\end{equation}
where $\tilde W \geq 0$ satisfies 
\begin{equation}\label{growth}
k_1\norm{\ts{F}}^p\le\tilde W( \mathbf{F}) \leq k_2[ \norm{\mathbf{F}}^p+1] \ \ 
\mathrm{for \ }
\mathbf{F}\in M^{m\times m}_+, \ p\in(m-1,m),
\end{equation}
for some positive constants $k_1,\,k_2$, and $h(\cdot)$ is a $C^2(0,\infty)$ 
convex function such that 
\begin{equation}\label{h}
h(\delta )\rightarrow \infty \mathrm{\ 
as \ }\delta \rightarrow 0^+, \quad \frac{h(\delta )}{\delta}
\rightarrow \infty \mathrm{\ as\ }\delta \rightarrow \infty.
\end{equation}
These hypotheses are typically satisfied by many stored energy functions which 
exhibit cavitating minimisers, for example,
\begin{equation}\label{SEtp}
W(\mathbf{F})=\mu\norm{\mathbf{F}}^p +h(\det\mathbf{F}), \quad \mu>0,
\end{equation}
where $h$ satisfies \eqref{h}.
\begin{remark}
It is well known that, under a variety of hypotheses (see, e.g, \cite{SiSp2})
on the stored energy function, there exists a minimiser of the energy (given by
\eqref{energy})
on the admissible set $\mathcal{A}_{\bx_0}$. Moreover, it is also known that if
$\mathbf{A}$ is sufficiently large\footnote{E.g., $\mathbf{A}=t\mathbf{B}$ for
some $\mathbf{B}\in M^{m\times m}_+$  with $t>0$ sufficiently large.}, then
any minimiser $\bu_0\in \mathcal{A}_{\bx_0}$ must satisfy $\alpha_{\bu_0}>0$ 
(see\cite{IUTAM}).
\end{remark}


We next prove that if we attempt to approximate, even in a weak sense, a 
singular deformation $\bu_0\in \mathcal{A}_{\bx_0}$ with finite elastic energy 
$E$ (given by \eqref{energy}) by a sequence of
nonsingular deformations in $\mathcal{A}_{\bx_0}$, then the energy of the 
approximating  sequence must necessarily diverge to infinity. In particular, 
this must also
hold in the case of approximating a singular energy minimiser. This phenomenon  
of the energy diverging to infinity is essentially due to the presence of the 
compressibility term $h$ which appears in the stored energy function 
\eqref{SEF}.

\begin{theorem}[The Repulsion Property]\label{thm:repprop}
Let $p\in (m-1, m)$. Suppose, for some $\mathbf{A}\in M^{m\times m}_+$,  that 
$\bu_0\in \mathcal{A}_{\bx_0}$ is a deformation with finite energy
and with $\alpha_{\bu_0}>0$. Suppose further that  $(\bu_n) \subset 
\mathcal{A}_{\bx_0}$ satisfies $\alpha_{\bu_n}=0, \ \forall n$ and that
$\bu_n\rightharpoonup \bu_0$ as $n\rightarrow \infty$ in $W^{1,p}(\Omega)$.  
Then $E(\bu_n)\rightarrow \infty$ as $n\rightarrow \infty$.
\end{theorem}
\begin{proof}
We first note that, since $||\mathbf{u}_n||< \mathrm{const.}$ uniformly  in $n$, 
it follows by \eqref{growth} that
\begin{equation*}
\mathrm{constant}\geq \int_\Omega \tilde W(\nabla \bu_n)\,\dif\ts{x} \ \ 
\mathrm{uniformly \ 
in \ } n.
\end{equation*}
We next claim that for any $R>0$ such that $B_R(\bx_0) \subset \Omega$ we have
\begin{equation*}
\int _{B_{R}(\bx_0)}\det (\nabla \bu_n)\,\dif\ts{x} \rightarrow \int 
_{B_{R}(\bx_0)}\det 
(\nabla \bu_0)\,\dif\ts{x}+\alpha _{\bu_0}>0, \mathrm{\ \ as\ \ }n\rightarrow 
\infty.
\end{equation*}
This follows from the facts (see \cite[Lemma 8.1]{MuSp}) that
\begin{eqnarray*}
(\mathrm{Det}(\nabla\bu_0))(B_R(\bx_0))& = &\int_{B_R(\bx_0)}\det (\nabla 
\bu_0) \,\dif\ts{x}
+\alpha_{\bu_0},\\
(\mathrm{Det}(\nabla\bu_n))(B_R(\bx_0))& = &\int_{B_R(\bx_0)}\det (\nabla
\bu_n)\,\dif\ts{x},\quad\mathrm{\ for \ all \ }n,
\end{eqnarray*}
and that
\begin{equation*}
(\mathrm{Det}(\nabla\bu_n))(B_R(\bx_0))\rightarrow
(\mathrm{Det}(\nabla\bu_0))(B_R(\bx_0))\mathrm{\ \ as\ \ }n\rightarrow \infty.
\end{equation*}
This last limit follows from classical results on the sequential weak 
continuity of the
mapping $\bu \rightarrow \mathrm{adj}(\nabla \bu)$ from $W^{1,p}$ into
$L^{\frac{p}{m-1}}$ (see \cite[Corollary 3.5]{Ba77}) and the compact embedding
of $W^{1,p}$ into $L^q_{\mbox{loc}}$ for every $q\in[1,\infty)$ for
functions satisfying the (INV) condition (see \cite[Lemma 3.3]{SiSp1}).

Hence, by Jensen's Inequality, for all $n$ we have
\begin{equation*}
E(\bu_n) \geq \int_{B_R(\bx_0)} W(\nabla \bu_n)\,\dif\ts{x} \geq
|B_R(\bx_0)|h\left(\frac{\int_{B_R(\bx_0)}\det(\nabla\bu_n\,\dif\ts{x})}
{|B_R(\bx_0)|}
\right) .
\end{equation*}
Hence
\begin{eqnarray*}
\mathrm{liminf}_{n\rightarrow\infty} E(\bu_n)\,\dif\ts{x} &\geq& 
\lim_{n\rightarrow \infty}
|B_R(\bx_0)|h\left(\frac{\int_{B_R(\bx_0)}\det(\nabla\bu_n)\,\dif\ts{x}}{
|B_R(\bx_0)|}
\right)\\
&=&|B_R(\bx_0)|h\left(\frac{\int_{B_R}\det(\nabla\bu_0)\,\dif\ts{x}+\alpha_{
\bu_0}}{|B_R|}
\right).
\end{eqnarray*}
Since this holds for all $R>0$ sufficiently small, and since $\alpha_{\bu_0}>0$ 
by assumption, it follows by \eqref{h} that
\begin{equation*}
\mathrm{liminf}_{n\rightarrow\infty} E(\bu_n)=\infty .
\end{equation*}

\end{proof}
\begin{remark}
If we replace the mode of convergence in the hypotheses of the above Theorem 
from weak convergence in $W^{1,p}$ to strong convergence, then it follows by the 
dominated convergence theorem that
\begin{equation*}
 \int_\Omega \tilde W(\nabla \bu_n)\,\dif\ts{x} \rightarrow \int_\Omega \tilde 
W(\nabla 
\bu_0)\,\dif\ts{x} \mathrm{\ \ as\ \ }n\rightarrow \infty.
\end{equation*}
Hence, this part of the total energy can be well approximated by nonsingular
deformations but the compressibility term involving $h$ cannot.
\end{remark}

\section{A decoupled method to circumvent the repulsion property}\label{sec:num}
We now consider an approximation scheme that avoids or works around
the repulsion property. The idea is to introduce a
decoupling or phase function $v$ in such a way that the difference
between the determinant of the approximation and the phase function
remains well behaved. The modified functional includes as well a
penalization term on $v$ reminiscent of the one used in the theory
of phase transitions, that penalizes if the function $v$ is not too large or 
not too small.

Let the stored energy function be as in \eqref{SEF}. For any $\tau>0$, let
$\phi_\tau:\Real\To[0,\infty)$ be a continuous function such that 
$\supp{\phi_\tau}\subset(0,\tau)$. For $\eps>0$, we define now the modified 
functional:
\begin{eqnarray}
I_\eps^\tau(\ts{u},v)&=&\int_{\Omega}\left[\tilde
W(\nabla\ts{u}(\ts{x}))+h(\det\nabla\ts{u}
(\ts{x})-v(\ts{x}))\right]\dif\ts{x}\nonumber\\
&&+\int_{\Omega}\left[\frac{\eps^\alpha}{\alpha} \norm{\nabla
v(\ts{x})}^\alpha+\frac{1}{q\eps^q}
\phi_\tau(v(\ts{x}))\right]\dif\ts{x},\label{modfunct}
\end{eqnarray}
where $\alpha>1$, $\frac{1}{\alpha}+\frac{1}{q}=1$, and 
$(\ts{u},v)\in\mc{U}$ where
\begin{eqnarray} \label{set:reg}
\mathcal{U} =\{(\bu,v)\in W^{1,p}(\Omega)\times W^{1,\alpha}(\Omega)
:\bu|_{\partial\Omega}=\bu^h,\
\bu_e\ {\rm satisfies}\ {\rm (INV)\ on}\ \Omega, \\ \nonumber
\det\nabla\bu>v\ge0\ {\rm a.e.},\ \
\mathrm{Det}\nabla\bu=(\det\nabla\bu)\mathcal{L}^m, \
v|_{\partial\Omega}=0\}.
\end{eqnarray}
The coupled $h$ term in this functional, because of \eqref{h}, penalizes for 
large $\det\nabla\bu$ and $v$ small. The term on $\nabla v$, for 
$\eps$ small, allows for large phase transitions in the function $v$. On the 
other hand, the term with the function $\phi_\tau$ for $\eps$ small, forces 
the regions where $v$ is positive but less than $\tau$, to have small 
measure, i.e. to ``concentrate''.

We now show that for any given $\tau,\eps>0$, the functional \eqref{modfunct} 
has a minimizer over $\mathcal{U}$.
\begin{lemma}\label{lem:epsF}
 Assume that $\tilde{W}(\cdot)$ and $h(\cdot)$ are nonnegative and that 
\eqref{growth}, \eqref{h} hold. For each $\tau>0$ and $\eps>0$ there 
exists 
$(\bu_\eps^\tau,v_\eps^\tau)\in\mc{U}$ such that
\[
I_\eps^\tau(\ts{u}_{\eps}^\tau,v_{\eps}^\tau)=\inf_{\mathcal{U}}I_\eps^\tau(\ts{
u },v).
\]
\end{lemma}
\begin{prueba}
 Since $\tilde{W}(\cdot)$ and $h(\cdot)$ are nonnegative and the pair 
$(\bu^h,0)$ belongs to $\mathcal{U}$, it
follows that $\inf_{\mathcal{U}}I_\eps^\tau(\ts{u},v)$ exists and (cf. 
\eqref{SEF})
\begin{equation}\label{aux:21}
 \inf_{\mathcal{U}}I_\eps^\tau(\ts{u},v)\le
I_\eps^\tau(\bu^h,0)= 
\int_{\Omega}W(\nabla\bu^h)\,\dif\ts{x}\equiv\ell.
\end{equation}
Let $\set{(\ts{u}_k,v_k)}$ be an infimizing sequence. From the above 
inequality, we can assume that $I_\eps^\tau(\bu_k,v_k)\le\ell$ for all $k$. It 
follows that 
\[
\int_{\Omega}\tilde{W}(\nabla\ts{u}_k(\ts{x}))\,\dif\ts{x}\le\ell,
\quad\forall k,
\]
which together with \eqref{growth} implies that for a subsequence 
$\set{\bu_k}$ 
(not relabelled), $\bu_k\rightharpoonup\bu_\eps^\tau$ in $W^{1,p}(\Omega)$, 
with $\bu_\eps^\tau=\bu^h$ over $\partial\Omega$ and $\bu_\eps^\tau$ satisfying 
the (INV) condition on $\Omega$.

From \eqref{aux:21} we get as well that
\[
\int_{\Omega}h(\det\nabla\ts{u}_k(\ts{x})-v_k(\ts{x}))\,\dif\ts{x}\le\ell,
\quad\forall k.
\]
This together with \eqref{h} and de la Vall\'ee Poussin criteria, imply that 
for a subsequence (not relabelled), $\det\nabla\ts{u}_k-v_k\rightharpoonup 
w_\eps^\tau$ in $L^1(\Omega)$, with $w_\eps^\tau>0$ a.e. Once again, 
\eqref{aux:21} implies (since $\eps$ is fixed) that $\set{v_k}$ is bounded in 
$W^{1,\alpha}(\Omega)$, and thus for a subsequence (not relabelled) that 
$v_k\rightharpoonup v_\eps^\tau$ in $W^{1,\alpha}(\Omega)$, with 
$v_\eps^\tau\ge0$ a.e. and $v_\eps^\tau=0$ on $\partial\Omega$. Thus we can 
conclude that $\det\nabla\ts{u}_k\rightharpoonup w_\eps^\tau+v_\eps^\tau$ in 
$L^1(\Omega)$. Since 
$\mathrm{Det}\nabla\ts{u}_{k}=(\det\nabla\ts{u}_{k})\mathcal{L}^m$, we
have that (see \cite[proof of Lemma (4.5)]{SiSpTi2006})
\[
 \wsconv{\det\nabla\ts{u}_{k}}{\mathrm{Det}\nabla\ts{u}_\eps^\tau}\quad\mbox{in 
}\Omega,
\]
from which it follows that 
$\mathrm{Det}\nabla\ts{u}_\eps^\tau=(w_\eps^\tau+v_\eps^\tau)\mc{L}^m$. Since 
$w_\eps^\tau+v_\eps^\tau\in L^1(\Omega)$, we have from \cite[Theorem 1]{Mu_A} 
that
\[ 
\mathrm{Det}\nabla\ts{u}_\eps^\tau=(\det\nabla\bu_\eps^\tau)\mc{L}^m,
\quad\det\nabla\bu_\eps^\tau=w_\eps^\tau+v_\eps^\tau.
\]
Thus $(\bu_\eps^\tau,v_\eps^\tau)\in\mc{U}$. Finally, since
\begin{eqnarray*}
& \bu_k\rightharpoonup\bu_\eps^\tau\mbox{ in }W^{1,p}(\Omega),\quad 
v_k\rightharpoonup v_\eps^\tau\mbox{ in }W^{1,\alpha}(\Omega),&\\ 
&\det\nabla\ts{u}_k-v_k\rightharpoonup 
w_\eps^\tau=\det\nabla\bu_\eps^\tau-v_\eps^\tau\mbox{ in 
}L^1(\Omega),&
\end{eqnarray*}
we have by the sequential weak lower semi--continuity of $I_\eps^\tau$, that
\[ 
I_\eps^\tau(\bu_\eps^\tau,v_\eps^\tau)\le\varliminf_{k\To\infty}I_\eps^\tau(\ts{
u}_{k },v_{ k} )= 
\inf_{\mathcal{U}}I_\eps^\tau(\ts{u},v),
\]
and thus
\[
 I_\eps^\tau(\bu_\eps^\tau,v_\eps^\tau)=\inf_{\mathcal{U}}I_\eps^\tau(\ts{u},v).
\]
\end{prueba}
Our next result shows that if $\ts{A}$ in \eqref{boundary}, is not too 
``large'', then the minimizer $(\bu_\eps^\tau,v_\eps^\tau)$ of Lemma 
\ref{lem:epsF} must be $(\bu^h,0)$.
\begin{proposition}\label{affine_global_gen}
Assume that the function $\tilde{W}$ is quasi--convex. If $\ts{A}$ in 
\eqref{boundary} is such that $\deru{h}(\det\ts{A})\le0$, then the global 
minimizer $(\bu_\eps^\tau,v_\eps^\tau)$ of $I_\eps^\tau(\cdot,\cdot)$ over
$\mathcal{U}$ is given by $\bu=\bu^h$ and $v=0$ in $\Omega$.
\end{proposition}
\begin{prueba}
Note that for any $(\bu,v)\in\mc{U}$, we have
\[
I_\eps^\tau(\bu,v)\ge 
\int_{\Omega}\left[\tilde
W(\nabla\ts{u}(\ts{x}))+h(\det\nabla\ts{u}
(\ts{x})-v(\ts{x}))\right]\dif\ts{x}.
\]
Since $\mathrm{Det}\nabla\bu=(\det\nabla\bu)\mathcal{L}^m$ and $\tilde{W}$ 
is quasi-convex, we have that
\[
\int_{\Omega}\tilde W(\nabla\ts{u}(\ts{x}))\dif\ts{x}\ge
\int_{\Omega}\tilde W(\nabla\ts{u}^h(\ts{x}))\dif\ts{x}.
\]
In addition, by the convexity of $h(\cdot)$ we get:
\[
h\left(\det\nabla\ts{u}(\ts{x})-v(\ts{x})\right)\ge 
h(\det\ts{A})+\deru{h}(\det\ts{A}) \left(
\det\nabla\ts{u}(\ts{x})-v(\ts{x})-\det\ts{A}\right).
\]
Hence
\begin{eqnarray*}
\int_{\Omega}h\left(\det\nabla\ts{u}(\ts{x})-v(\ts{x})\right)\dif\ts{x}
&\ge& 
\int_{\Omega}h(\det\ts{A})\dif\ts{x}
-\deru{h}(\det\ts{A})\int_{\Omega}v(\ts{x})\dif\ts{x}\\
&&+\deru{h}(\det\ts{A})\int_{\Omega}\left(\det\nabla\ts{u}(\ts{x})
-\det\ts{A}\right)\dif\ts{x}
\end{eqnarray*}
Again, since $\mathrm{Det}\nabla\bu=(\det\nabla\bu)\mathcal{L}^m$, we have that
\[
\int_{\Omega}\left(\det\nabla\ts{u}(\ts{x})
-\det\ts{A}\right)\dif\ts{x}=0.
\]
Using now that $\deru{h}(\det\ts{A})\le0$ and that $v\ge0$, we get
\[
\int_{\Omega}h\left(\det\nabla\ts{u}(\ts{x})-v(\ts{x})\right)\dif\ts{x}\ge 
\int_{\Omega}h(\det\ts{A})\dif\ts{x}.
\]
Combining this with the two inequalities at the beginning of this proof, we get 
that
\begin{eqnarray*}
I_\eps^\tau(\bu,v)&\ge&\int_{\Omega}\left[\tilde 
W(\nabla\ts{u}^h(\ts{x}))+h(\det\nabla\ts{u}^h(\ts{x}))\right]\dif\ts{x}
=I_\eps^\tau(\bu^h,0).
\end{eqnarray*}
Since $(\bu,v)$ is arbitrary in $\mc{U}$ and 
$(\bu^h,0)\in\mc{U}$, we have that $(\bu^h,0)$ is the global
minimizer in this case.
\end{prueba}

Let $\mc{M}(\Omega)$ be the space of signed Radon measures on $\Omega$. If 
$\mu\in\mc{M}(\Omega)$, then
\[
 \langle\mu,v \rangle=\int_\Omega v\,\dif\mu,\quad\forall\,v\in C_0(\Omega),
\]
where $C_0(\Omega)$ denotes the set of continuous functions with compact 
support in $\Omega$. Moreover
\[
\norm{\mu}_{\mc{M}(\Omega)}=\sup\set{|{\langle\mu,v \rangle}|\,:\,v\in 
C_0(\Omega),\,\norm{v}_{L^\infty(\Omega)}\le1}.
\]
A sequence $\set{\mu_n}$ in $\mc{M}(\Omega)$ converges weakly $*$ to 
$\mu\in\mc{M}(\Omega)$, denoted $\mu_n\stackrel{*}{\rightharpoonup}\mu$, if
\[
 \limite{n}{\infty}\langle\mu_n,v \rangle= \langle\mu,v \rangle,\quad\forall 
\,v\in C_0(\Omega).
\]
Note that any function in $L^1(\Omega)$ can be regarded as belonging to 
$\mc{M}(\Omega)$ with the same norm. It follows from this observation and the 
weak compactness of $\mc{M}(\Omega)$, that if $\set{v_n}$ is a bounded sequence 
in $L^1(\Omega)$, then it has a subsequence $\set{v_{n_k}}$ such that 
$v_{n_k}\stackrel{*}{\rightharpoonup}\mu$ where $\mu\in\mc{M}(\Omega)$.

For any subset $E$ of $\Omega$, we define its \textit{(Caccioppoli) perimeter 
in $\Omega$} by
\[ 
P(E,\Omega)=\sup\left\{\int_{\Omega}\chi_E(\ts{x})\,\mathrm{div}\,
\tsg{\phi}(\ts{x})\,\dif\ts{x}\,:\,\tsg{\phi}\in C^1(\Omega;\Real^m),\,
\norm{\tsg{\phi}}_{L^\infty(\Omega)}\le1\right\}.
\]
$E$ is said to have \textit{finite perimeter} in $\Omega$ if 
$P(E,\Omega)<\infty$. For a set of finite perimeter, it follows from the 
Gauss--Green Theorem (cf. \cite[Thm. 5.16]{EvGa2015}) that
\[
 P(E,\Omega)=\mc{H}^{m-1}(\partial_*E),
\]
where $\partial_*E$ is the so called \textit{measure theoretic boundary} of $E$.

We now study the convergence of the minimizers in Lemma \ref{lem:epsF} as 
$\eps\To0$. We employ the following notation:
\[
 H_\tau(s)=\int_0^s\phi_\tau^{1/q}(t)\,\dif t.
\]
Using this we can now prove the following:
\begin{theorem}\label{DPconvG}
 Assume a stored energy of the form \eqref{SEF}--\eqref{h} and that 
$p\in(m-1,m)$. Let 
$(\bu_\eps^\tau,v_\eps^\tau)\in\mc{U}$ be a minimizer of $I_\eps^\tau$ over 
$\mc{U}$. Then for any sequence $\eps_j\To0$, the sequences 
$\set{\ts{u}_j^\tau}$ and $\set{v_j^\tau}$, where 
$\ts{u}_j^\tau=\ts{u}_{\eps_j}^\tau$ and 
$v_j^\tau=v_{\eps_j}^\tau$, have subsequences $\set{\ts{u}_{j_k}^\tau}$ and 
$\set{v_{j_k}^\tau}$ with
$\ts{u}_{j_k}^\tau\rightharpoonup\ts{u}^\tau$ in $W^{1,p}(\Omega)$ and
$\wsconv{v_{j_k}^\tau}{\nu^\tau}$ in $\mc{M}(\Omega)$ where $\nu^\tau$ is a 
nonnegative Radon measure. Moreover
$\bu^\tau|_{\partial\Omega}=\bu^h$, $\bu_e^\tau$ satisfies (INV) on $\Omega$, 
and
\[
\mathrm{Det}\nabla\bu^\tau=(\det\nabla\bu^\tau)\mc{L}^m+\nu_s^\tau,
\]
where $\det\nabla\bu^\tau\in L^1(\Omega)$ with $\det\nabla\bu^\tau>0$ a.e. in 
$\Omega$ and $\nu_s^\tau$ is the singular part of $\nu^\tau$ with respect to 
Lebesgue measure. If we let 
$\hat{v}_{j_k}^\tau(\ts{x})=\min\set{v_{j_k}^\tau(\ts{x}),\tau}$, then
$\set{\hat{v}_{j_k}^\tau}$ has a subsequence converging in $L^1(\Omega)$ to a
function $g^\tau$ that assumes only the values $0$ and $\tau$ a.e., and
\begin{equation}\label{ineqG}
\varliminf_{k\To\infty}I_{\eps_{j_k}}^\tau(\ts{u}_{j_k}^\tau,v_{j_k}^\tau)\ge 
 \int_{\Omega}\left[\tilde
W(\nabla\ts{u}^\tau(\ts{x}))+h(\det\nabla\ts{u}^\tau
(\ts{x})-\omega^\tau(\ts{x}))\right]\dif\ts{x}
+H_\tau(\tau)P(B_\tau,\Omega),
\end{equation}
where $\omega^\tau\in L^1(\Omega)$ is the derivative of $\nu^\tau$ with respect 
to Lebesgue measure and satisfies that $\det\nabla\ts{u}^\tau>\omega^\tau\ge0$ 
a.e.,  and  $B_\tau=\set{\ts{x}\in\Omega\,:\,g^\tau(\ts{x})=0}$.
\end{theorem}
\begin{prueba}
 The inequality
 \begin{equation}\label{unif_bdd}
I_{\eps_j}^\tau(\ts{u}_j^\tau,v_j^\tau)\le\int_{\Omega}W(\nabla\bu^h)\,\dif\ts{x
},
\end{equation}
together with \eqref{growth} and Poincar\'e's inequality, imply the existence 
of 
a subsequence $\set{\ts{u}_{j_k}^\tau}$ converging weakly to a function 
$\ts{u}^\tau$ in $W^{1,p}(\Omega)$. Clearly $\bu^\tau|_{\partial\Omega}=\bu^h$, 
and that $\bu_e^\tau$ satisfies (INV) on $\Omega$ follows from \cite[Lemma 
3.3]{MuSp}. From \eqref{h} and de la Vall\'ee Poussin criteria, it follows that 
there is a subsequence (with indexes written as for the previous one) 
$\set{\det\nabla\ts{u}_{j_k}^\tau-v_{j_k}^\tau}$ such that
\begin{equation}\label{convgD}
 \det\nabla\ts{u}_{j_k}^\tau-v_{j_k}^\tau\rightharpoonup w^\tau,\quad\mbox{in 
}L^1(\Omega).
\end{equation}
Since $\det\nabla\ts{u}_{j_k}^\tau-v_{j_k}^\tau>0$ a.e. on $\Omega$, the first 
condition in \eqref{h} implies that we must have that $w^\tau>0$ a.e. on 
$\Omega$. Now from $\det\nabla\ts{u}_{j_k}^\tau>v_{j_k}^\tau\ge0$, it follows 
that
\[
 \int_\Omega v_{j_k}^\tau\,\dif\ts{x}\le\int_\Omega 
\det\nabla\ts{u}_{j_k}^\tau\,\dif\ts{x}=|\ts{u}^h(\Omega)|.
\]
Thus $\set{v_{j_k}^\tau}$ is bounded in $L^1(\Omega)$. Hence there exists 
$\nu^\tau\in\mc{M}(\Omega)$ such that (for a subsequence denoted the same) 
$\wsconv{v_{j_k}^\tau}{\nu^\tau}$ in $\mc{M}(\Omega)$. Since $v_{j_k}^\tau\ge0$ 
for all $k$, the measure $\nu^\tau$ must be non--negative. Combining this with 
\eqref{convgD} we get that
\[
\wsconv{(\det\nabla\ts{u}_{j_k}^\tau)\mc{L}^m}{w^\tau\mc{L}^m+\nu^\tau}
\quad\mbox{in }\Omega.
\]
Since 
$\mathrm{Det}\nabla\ts{u}_{j_k}^\tau=(\det\nabla\ts{u}_{j_k}^\tau)\mathcal{L}^m$
, we have that (see \cite[proof of Lemma 4.5]{SiSpTi2006})
\[ 
\wsconv{(\det\nabla\ts{u}_{j_k}^\tau)\mc{L}^m}{\mathrm{Det}\nabla\ts{u}^\tau}
\quad\mbox{in }\Omega,
\]
from which it follows that 
$\mathrm{Det}\nabla\ts{u}^\tau=w^\tau\mc{L}^m+\nu^\tau$. By the
Lebesgue decomposition theorem, $\nu^\tau=\nu_{ac}^\tau+\nu_s^\tau$ where 
$\nu_{ac}^\tau$ is absolutely continuous with respect to $\mc{L}^m$ and 
$\nu_s^\tau$ is singular with respect to $\mc{L}^m$. Thus 
$\mathrm{Det}\nabla\ts{u}^\tau=w^\tau\mc{L}^m+\nu_{ac}^\tau+\nu_s^\tau$. 
Since $w^\tau\mc{L}^m+\nu_{ac}^\tau$ is absolutely continuous with respect to 
$\mc{L}^m$, it follows by the uniqueness in the Lebesgue decomposition theorem, 
that $w^\tau\mc{L}^m+\nu_{ac}^\tau$ is the absolutely continuous part of 
$\mathrm{Det}\nabla\ts{u}^\tau$ with respect to $\mc{L}^m$. Since $p>m-1$
and $\bu_e^\tau$ satisfies (INV) on $\Omega$, the conclusions of Theorem 1 in
\cite{Mu_A} hold. In particular, from Remark 2 of that theorem, we get that the 
absolutely continuous part of $\mathrm{Det}\nabla\ts{u}^\tau$ is 
$(\det\nabla\ts{u}^\tau)\mc{L}^m$. Thus, by the uniqueness in the Lebesgue 
decomposition theorem, we must have that 
$(\det\nabla\ts{u}^\tau)\mc{L}^m=w^\tau\mc{L}^m+\nu_{ac}^\tau$. Hence 
$\mathrm{Det}\nabla\ts{u}^\tau=(\det\nabla\ts{u}^\tau)\mc{L}^m+\nu_s^\tau$ and 
$\det\nabla\ts{u}^\tau=w^\tau+\omega^\tau$ where $\omega^\tau$ is the 
derivative of $\nu^\tau$ with respect to $\mc{L}^m$. Since $w^\tau>0$ and 
$\omega^\tau\ge0$ a.e., it follows that $\det\nabla\ts{u}^\tau 
>\omega^\tau\ge0$ a.e.

Since $\phi_\tau$ is nonnegative and $\supp{\phi_\tau}\subset(0,\tau)$, it 
follows that $\set{H_\tau(v_{j_k}^\tau)}$ is bounded in $L^1(\Omega)$. Moreover
\[
 \int_{\Omega}\left[\frac{\eps_{j_k}^\alpha}{\alpha} \norm{\nabla
v_{j_k}^\tau(\ts{x})}^\alpha+\frac{1}{q\eps_{j_k}^q}\phi_\tau(v_{j_k}^\tau 
(\ts{x}))\right]\dif\ts{x}\ge 
\int_\Omega\norm{\nabla[ H_\tau(v_{j_k}^\tau(\ts{x}))]}\,\dif\ts{x}.
\]
If we let $\hat{v}_{j_k}^\tau(\ts{x})=\min\set{v_{j_k}^\tau(\ts{x}),\tau}$, then
\[
 \int_\Omega\norm{\nabla[ H_\tau(v_{j_k}^\tau(\ts{x}))]}\,\dif\ts{x}= 
\int_\Omega\norm{\nabla[ H_\tau(\hat{v}_{j_k}^\tau(\ts{x}))]}\,\dif\ts{x}
\]
It follows $\set{H_\tau(\hat{v}_{j_k}^\tau)}$ is bounded in $BV(\Omega)$ (cf. 
\cite{Mo1987}) and thus it has a subsequence converging in $L^1(\Omega)$. Since 
$\hat{v}_{j_k}^\tau:\Omega\To[0,\tau]$, we get that $\hat{v}_{j_k}^\tau\To 
g^\tau$ in $L^1(\Omega)$. In addition 
\begin{eqnarray*}
 \int_\Omega\phi_\tau(v_{j_k}^\tau(\ts{x}))\,\dif\ts{x}&=& 
\int_\Omega\phi_\tau(\hat{v}_{j_k}^\tau(\ts{x}))\,\dif\ts{x},\\ 
\limite{k}{\infty}\int_\Omega\phi_\tau(v_{j_k}^\tau(\ts{x}))\,\dif\ts{x}&=&0,
\quad\mbox{(cf. \eqref{unif_bdd})},
\end{eqnarray*}
from which we get that $\int_\Omega\phi_\tau(g^\tau(\ts{x}))\,\dif\ts{x}=0$, 
i.e., that $g^\tau$ assumes only the values $0$ or $\tau$ a.e. Also
\begin{multline*}
 \varliminf_{k\To\infty}\int_{\Omega}\left[\frac{\eps_{j_k}^\alpha}{
\alpha}\norm{\nabla v_{j_k}^\tau(\ts{x})}^\alpha
+\frac{1}{q\eps_{j_k}^q}\phi_\tau(v_{j_k}^\tau(\ts{x}))\right]\dif\ts{x}
\ge\\
\varliminf_{k\To\infty}\int_\Omega\norm{\nabla[ 
H_\tau(\hat{v}_{j_k}^\tau(\ts{x}))]}\,\dif\ts{x}\ge 
\int_\Omega\norm{\nabla[H_\tau(g^\tau(\ts{x}))]}\,\dif\ts{x}
=H_\tau(\tau)P(B_\tau,\Omega),
\end{multline*}
where for the second inequality we used the lower semicontinuity property of 
the variation measure (cf. \cite[Thm. 5.2]{EvGa2015}).
Finally combining this result with those from the first part 
of this proof and the weak lower semicontinuity property of the mechanical part 
of the functional \eqref{modfunct}, we get that \eqref{ineqG} follows.
\end{prueba}

Note that Theorem \ref{DPconvG} in a sense falls short of fully characterizing 
any possible singular behaviour in a minimizer $\ts{u}^*$ of the energy 
functional \eqref{energy}. Since the parameter $\tau$ is fixed, the phase 
functions are not ``forced'' to follow or mimic the singular behaviour in 
$\ts{u}^*$ once they have crossed the barrier $\tau$. Moreover, the actual 
location of the set of possible singularities  in $\ts{u}^*$ has not been fully 
resolved due to the presence of the function $\omega^\tau$ in the $h$--term of 
the energy functional. Thus we need to study the behaviour of the functions 
$\ts{u}^\tau$, $\omega^\tau$, $g^\tau$, and the 
measures $\nu^\tau$ as $\tau\To\infty$. 

In the sequel we employ some of the 
notation within the proof of 
Theorem \ref{DPconvG} as well as the following: given $\tau_1>0$ and a 
sequence $\set{\eps_j}$ converging to zero, we apply Theorem \ref{DPconvG} to 
get a subsequence $\set{\eps_{1,r}}$ of $\set{\eps_j}$ with the corresponding 
sequences of functions $\set{\ts{u}_{1,r}}$, $\set{v_{1.r}}$, etc. We keep 
denoting the limiting functions and measures by $\ts{u}^{\tau_1}$, 
$\nu^{\tau_1}$, etc. Now given any $\tau_k$ with $k>1$, we apply Theorem 
\ref{DPconvG} to the subsequence $\set{\eps_{k-1,r}}$ obtained from 
$\tau_{k-1}$, to get a new subsequence $\set{\eps_{k,r}}$ of 
$\set{\eps_{k-1,r}}$, and so on. After relabelling, we denote by 
$\set{\ts{u}_{k,r}}$, $\set{v_{k,r}}$, etc., the sequences obtained from 
Theorem \ref{DPconvG} by this process for any given $\tau_k$.

\begin{lemma}\label{lem:4.3a}
 The sequences $\set{g^{\tau_k}}$ ans $\set{\nu^{\tau_k}}$ have subsequences 
(relabelled the same) such for some $\nu,\nu^*\in\mc{M}(\Omega)$, we have 
$\wsconv{g^{\tau_k}}{\nu}$ and $\wsconv{\nu^{\tau_k}}{\nu^*}$ in 
$\mc{M}(\Omega)$.
\end{lemma}
\begin{prueba}
Note that
\[
 \int_\Omega \hat{v}_{k,r}(\ts{x})\,\dif\ts{x}\le \int_\Omega 
v_{k,r}(\ts{x})\,\dif\ts{x}\le|\ts{u}^h(\Omega)|,
\]
and since $\hat{v}_{k,r}^\tau\To g^{\tau_k}$ in $L^1(\Omega)$ as $r\To\infty$, 
it follows that
\[
 \int_\Omega 
g^{\tau_k}(\ts{x})\,\dif\ts{x}\le|\ts{u}^h(\Omega)|,\quad\forall k.
\]
Thus for some subsequence of $\set{\tau_k}$ (not relabelled), we have 
that $\wsconv{g^{\tau_k}}{\nu}$ in $\mc{M}(\Omega)$, for some 
$\nu\in\mc{M}(\Omega)$. 

Also, since $\wsconv{v_{k,r}}{\nu^{\tau_k}}$ as $r\To\infty$, we get 
that for any $\phi\in C_0(\Omega)$, $\norm{\phi}_{L^\infty(\Omega)}\le1$, we 
have that
\[
 \limite{r}{\infty}\int_\Omega v_{k,r}(\ts{x})\phi(\ts{x})\,\dif\ts{x}
= \langle\nu^{\tau_k},\phi \rangle.
\]
But
\[
\left|\int_\Omega 
v_{k,r}(\ts{x})\phi(\ts{x})\,\dif\ts{x}\right|\le\int_\Omega 
v_{k,r}(\ts{x})\,\dif\ts{x}\le|\ts{u}^h(\Omega)|.
\]
Letting $r\To\infty$ we get that $\abs{\langle\nu^{\tau_k},\phi \rangle}\le 
|\ts{u}^h(\Omega)|$, and hence that $\norm{\nu^{\tau_k}}_{\mc{M}(\Omega)}\le 
|\ts{u}^h(\Omega)|$. Thus by taking a subsequence of $\set{\tau_k}$ 
(relabelled the same), we have $\wsconv{\nu^{\tau_k}}{\nu^*}$ in 
$\mc{M}(\Omega)$, 
for some $\nu^*\in\mc{M}(\Omega)$.
\end{prueba}
From these results and 
\cite[Thm. 5.1]{Conway1990}, we get the following:

\begin{lemma}\label{lem:4.3}
 The sequences $\set{\hat{v}_{k,r}}$ and $\set{v_{k,r}}$ have 
subsequences $\set{\hat{v}_k}$ and $\set{v_k}$ respectively, where 
$\hat{v}_k=\hat{v}_{k,r_k}$ and $v_k=v_{k,r_k}$ with $r_k\To\infty$, such that 
$\wsconv{\hat{v}_{k}}{\nu}$ and $\wsconv{v_{k}}{\nu^*}$ in 
$\mc{M}(\Omega)$, as $\tau_k\To\infty$.
\end{lemma}
We now show that $\nu=\nu^*$ and that $\nu$ is singular with respect to 
$\mc{L}^m$.
\begin{lemma}\label{lem:4.6}
If $\set{\tau_k}$ is such that $\tau_k\To\infty$, then 
$\wsconv{v_{k}-\hat{v}_k}{0}$ in $\mc{M}(\Omega)$ and consequently 
that $\nu=\nu^*$.
\end{lemma}
\begin{prueba}
 Let
 \[
  E_k=\set{\ts{x}\in\Omega\,:\,v_k(\ts{x})\ge\tau_{k}}.
 \]
Since
\[
 \int_{\Omega}v_k(\ts{x})\,\dif\ts{x}\ge 
\int_{E_n}v_k(\ts{x})\,\dif\ts{x}\ge\tau_{k}|E_k|,
\]
we get that
\[
 |E_k|\le\frac{1}{\tau_k}|\ts{u}^h(\Omega)|.
\]
Since $v_k-\hat{v}_k=(v_k-\tau_{k})\chi_{E_k}$, it follows that for any 
$\psi\in C_0(\Omega)$ with $\norm{\psi}_{L^\infty(\Omega)}\le1$:
\[
 \left|\int_\Omega\psi(\ts{x})(v_k(\ts{x})-\hat{v}_k(\ts{x}))\,\dif\ts{x} 
\right| 
\le\norm{\psi}_{L^\infty(\Omega)}\int_{E_k}v_k(\ts{x})\,\dif\ts{x}\To0,\quad 
\mbox{ as }k\To\infty,
\]
as $\int_{\Omega}v_k\,\dif\ts{x}\le|\ts{u}^h(\Omega)|$ and $|E_k|\To0$ as 
$k\To\infty$. Thus $\wsconv{v_{k}-\hat{v}_k}{0}$ in $\mc{M}(\Omega)$ and by 
Lemma \ref{lem:4.3} we have that $\nu=\nu^*$.
\end{prueba}

\begin{proposition}\label{prop:nusing}
There exist sets $B$ and $D$ disjoint such that $\Omega=B\cup D$, where $|D|=0$ 
and $\nu(B)=0$, i.e., $\nu$ is singular with respect to Lebesgue measure 
$\mc{L}^m$.
\end{proposition}
\begin{prueba}
 Let
 \[
  B_k=\set{\ts{x}\in\Omega\,:\,g^{\tau_{k}}(\ts{x})=0}.
 \]
Without loss of generality, we can assume that $\tau_k\ge k^2$. Since 
$g^{\tau_{k}}$ assumes only the values $0$ and $\tau_{k}$, we now have that
\[
 |\Omega\setminus B_k|\le\frac{1}{k^2}|\ts{u}^h(\Omega)|.
\]
Let
\[
 C_M=\bigcap_{j\ge M}^\infty B_j,\quad B=\bigcup_{M=1}^\infty C_M.
\]
Note that $C_M\subset C_{M+1}$ for all $M$. Also, let
\[
 D_M=\bigcup_{j\ge M}^\infty (\Omega\setminus B_j),\quad D=\bigcap_{M=1}^\infty 
D_M.
\]
Note that $D_M=\Omega\setminus C_M$ and that $D_{M+1}\subset D_M$ for all $M$. 
Also $D=\Omega\setminus B$. Since
\[
 |D_M|\le\sum_{j=M}^\infty|\Omega\setminus B_j|\le|\ts{u}^h(\Omega)|\sum_{j=M} 
^\infty\frac{1}{j^2},
\]
we have that
\[
 |D|=\limite{M}{\infty}|D_M|=0.
\]
Thus $|B|=|\Omega|$. We now show that $\nu(B)=0$. For any $M\ge1$, we have
\[
 \nu(B)=\nu(B\cap C_M)+\nu(B\setminus C_M).
\]
Since $\wsconv{g^{\tau_{k}}}{\nu}$ in $\mc{M}(\Omega)$, we have that
\[
 \nu(B\cap C_M)=\limite{n}{\infty}\int_{B\cap C_M}g^{\tau_{n}}(\ts{x}) 
\,\dif\ts{x}.
\]
Since $C_M\subset C_k$ for all $k>M$ and $C_k\subset B_k$, it follows that 
\[
 \int_{B\cap C_M}g^{\tau_{k}}(\ts{x}) 
\,\dif\ts{x}=0,\quad k>M.
\]
Hence  $\nu(B\cap C_M)=0$ for all $M\ge1$ and we get now that
\begin{eqnarray*}
 \nu(B)&=&\limite{M}{\infty}[\nu(B\cap C_M)+\nu(B\setminus C_M)]\\
 &=&\limite{M}{\infty}\nu(B\setminus C_M)
 =\limite{M}{\infty}\nu(B\cap D_M)\\
 &=&\nu(B\cap D)=0.
\end{eqnarray*}
Thus we have $\Omega=B\cup D$, $B$ and $D$ disjoint, with $|D|=0$ and 
$\nu(B)=0$, i.e. $\nu$ is singular with respect to Lebesgue measure $\mc{L}^m$.
\end{prueba}
We now give the corresponding results for the sequences $\set{\bu^{\tau_k}}$ 
and $\set{\omega^{\tau_k}}$.
\begin{proposition}\label{wconv_utau}
 Let $\set{\tau_k}$ be a sequence such that $\tau_k\To\infty$. Then the 
sequences $\set{\bu^{\tau_k}}$ and $\set{\omega^{\tau_k}}$ have subsequences 
relabelled the same, such that $\ts{u}^{\tau_k}\rightharpoonup\ts{u}^*$ in 
$W^{1,p}(\Omega)$, $\det\nabla\ts{u}^{\tau_{k}}\rightharpoonup 
\det\nabla\ts{u}^*$ in $L^1(\Omega)$, and $\omega^{\tau_k}\To0$ in 
$L^1(\Omega)$. Moreover, the function $\bu^*$ is such that 
$\ts{u}^*|_{\partial\Omega}=\bu^h$, $\bu_e^*$ satisfies 
(INV) on $\Omega$, and
\begin{equation}\label{DETlim}
\mathrm{Det}\nabla\ts{u}^*=(\det\nabla\ts{u}^*)\mc{L}^m+\nu,
\end{equation}
with $\det\nabla\ts{u}^*\in L^1(\Omega)$ and $\det\nabla\ts{u}^*>0$ a.e. in 
$\Omega$.
\end{proposition}
\begin{prueba}
  Since, by Lemma \ref{lem:4.3a}, $\wsconv{\nu^{\tau_{k}}}{\nu}$, we have that 
$\nu^{\tau_{k}}(B)\To\nu(B)=0$, where $B$ is as in Proposition 
\ref{prop:nusing}. But
\[
 \int_B\omega^{\tau_{k}}\,\dif\ts{x}\le \nu^{\tau_{k}}(B).
\]
As $\omega^{\tau_{k}}\ge0$ a.e., the above implies that
\[
 \limite{k}{\infty}\int_B\omega^{\tau_{k}}\,\dif\ts{x}=0,
\]
which implies that $\omega^{\tau_{k}}\To0$ in $L^1(\Omega)$, where we used 
that $|B|=|\Omega|$.

From \eqref{growth}, \eqref{ineqG}, \eqref{unif_bdd}, and Poincar\'e's 
inequality, we get that for a subsequence of $\set{\bu^{\tau_k}}$ (not 
relabelled), we have $\ts{u}^{\tau_k}\rightharpoonup\ts{u}^*$ in 
$W^{1,p}(\Omega)$ for some function $\bu^*\in W^{1,p}(\Omega)$. Clearly 
$\bu^*|_{\partial\Omega}=\bu^h$, 
and that $\bu_e^*$ satisfies (INV) on $\Omega$ follows from \cite[Lemma 
3.3]{MuSp} and the fact that each $\ts{u}_k$ satisfies (INV).

From \eqref{h} and de la Vall\'ee Poussin criteria, it follows that 
there is a subsequence (with indexes written as for the previous one) 
$\set{\det\nabla\ts{u}^{\tau_{k}}-\omega^{\tau_{k}}}$ such that
\[
 \det\nabla\ts{u}^{\tau_{k}}-\omega^{\tau_{k}}\rightharpoonup 
w^*,\quad\mbox{in }L^1(\Omega).
\]
Since $\det\nabla\ts{u}^{\tau_{k}}-\omega^{\tau_{k}}>0$ a.e. on $\Omega$, 
the first condition in \eqref{h} implies that we must have that $w^*>0$ a.e. 
on $\Omega$. Now $\det\nabla\ts{u}^{\tau_{k}}>\omega^{\tau_{k}}\ge0$ a.e. 
on $\Omega$, and since $\omega^{\tau_{k}}\To0$ in $L^1(\Omega)$, we get from 
the previous convergence that
 \[
 \det\nabla\ts{u}^{\tau_{k}}\rightharpoonup 
w^*,\quad\mbox{in }L^1(\Omega).
\]
It follows now from \cite[Theorem 4.2]{MuSp}, that $\det\nabla\ts{u}^*=w^*$. 
From the proof of Theorem \ref{DPconvG}, we have that 
$\det\nabla\ts{u}^{\tau_{k}}=w^{\tau_{k}}+\omega^{\tau_{k}}$ from which 
it follows that $w^{\tau_{k}}\rightharpoonup w^*$, in $L^1(\Omega)$. Also 
$\mathrm{Det}\nabla\ts{u}^{\tau_{k}}=w^{\tau_{k}}\mc{L}^m+\nu^{\tau_{k}}$ 
and since $\wsconv{\mathrm{Det}\nabla\ts{u}^{\tau_{k}}} 
{\mathrm{Det}\nabla\ts{u}^*}$, we get that \eqref{DETlim} holds.
\end{prueba}
We now have one of the main results of this paper.
\begin{theorem}\label{reppnumsch}
Let $\set{\tau_k}$ and $\set{\eps_r}$ be sequences such that 
$\tau_k\To\infty$ and $\eps_r\To0^+$, and let $(\bu_{k,r},v_{k,r})$ be a 
minimizer of $I_{\eps_{r}}^{\tau_k}$ over $\mc{U}$. Then there exist a 
subsequence of $\set{\tau_k}$ relabelled the same, and a subsequence 
$\set{\eps_{r_k}}$, such that if $(\bu_k,v_k)=(\bu_{k,r_k},v_{k,r_k})$, then  
$\ts{u}_k\rightharpoonup\ts{u}^*$ in $W^{1,p}(\Omega)$ and 
$\wsconv{v_{k}}{\nu}$ in $\mc{M}(\Omega)$ as 
$k\To\infty$. Moreover, if 
\begin{equation}\label{aux:p3}
\varliminf_{k\To\infty}H_{\tau_{k}}(\tau_{k})\ge c,
\end{equation}
for some constant $c>0$, then
\begin{equation}\label{ineqG2}
\varliminf_{k\To\infty}I_{\eps_{r_k}} 
^{\tau_{k}}(\ts{u}_{k},v_{k})\ge 
 \int_{\Omega}W(\nabla\ts{u}^*(\ts{x}))\,\dif\ts{x}+c\,P(B,\Omega).
\end{equation}
\end{theorem}
\begin{prueba}
The existence and the convergence of the subsequence $\set{v_k}$ with 
$v_k=v_{k,r_k}$, follows from the boundedness of $\set{v_{k,r}}$ in 
$L^1(\Omega)$, Theorem \ref{DPconvG}, Lemmas \ref{lem:4.3a} and \ref{lem:4.6}, 
and \cite[Thm. 5.1]{Conway1990}. For the existence and the convergence of the 
subsequence $\set{\bu_k}$ with $\bu_k=\bu_{k,r_k}$, it follows from the 
boundedness of this sequence in $W^{1,p}(\Omega)$ (cf. \eqref{unif_bdd}), 
Theorem \ref{DPconvG}, Proposition \ref{wconv_utau}, and \cite[Thm. 
5.1]{Conway1990}. 

Without loss of generality, we can assume that for each $k$, the $r_k$ is 
chosen so that
\[
 I_{\eps_{k,r_k}}^{\tau_{k}}(\ts{u}_{k,r_k},v_{k,r_k})>
\varliminf_{r\To\infty}I_{\eps_{k,r}}^{\tau_{k}}
(\ts{u}_{k,r},v_{k,r})-\frac{1}{k}.
\]
We get now using \eqref{ineqG} that
\begin{eqnarray}
I_{\eps_{k,r_k}}^{\tau_{k}}(\ts{u}_{k},v_{k})&\ge& 
 \int_{\Omega}\left[\tilde W(\nabla\ts{u}^{\tau_{k}}(\ts{x}))
+h(\det\nabla\ts{u}^{\tau_{k}}(\ts{x})
-\omega^{\tau_{k}}(\ts{x}))\right]\dif\ts{x}\nonumber\\
&&+H_{\tau_{k}}(\tau_{k})P(B_k,\Omega)-\frac{1}{k}.\label{aux:p1}
\end{eqnarray}
Since $B=\varliminf_kB_k$ it follows that $\chi_B=\varliminf_k\chi_{B_k}$. Thus 
for 
any $\tsg{\phi}\in C^1(\Omega;\Real^n)$, with 
$\norm{\tsg{\phi}}_{L^\infty(\Omega)}\le1$, we have that (cf. \cite[Ex. 12, 
Pag. 90]{Ro1968})
\[
 \int_\Omega \chi_B(\ts{x})\,\mathrm{div}\,\tsg{\phi}(\ts{x})\,\dif\ts{x}
 \le\varliminf_k\int_\Omega \chi_{B_k}(\ts{x})
 \,\mathrm{div}\,\tsg{\phi}(\ts{x})\,\dif\ts{x}\le
 \varliminf_k P(B_k,\Omega),
\]
from which we get that
\begin{equation}\label{aux:p2}
 P(B,\Omega)\le\varliminf_k P(B_k,\Omega).
\end{equation}
The result \eqref{ineqG2} now follows from \eqref{aux:p3}, \eqref{aux:p1}, 
\eqref{aux:p2}, and the convergence results in Proposition \ref{wconv_utau} for 
the sequences 
$\set{\ts{u}^{\tau_{k}}}$, $\set{\det\nabla\ts{u}^{\tau_{k}}}$, and
$\set{\omega^{\tau_{k}}}$.
\end{prueba}

The measure $\nu$ in this theorem, according to Proposition \ref{prop:nusing}, 
is concentrated on the set $D$ which is the complement of $B$. In addition, by 
the extended Lebesgue Decomposition Theorem (cf \cite{Ha1974}, 
\cite{HeSt1975}), $\nu$ is the sum of a discrete measure and a continuous one, 
both singular with respect to Lebesgue measure. The discrete part of $\nu$ 
correspond to points in the reference configuration where singularities of 
cavitation type may occur, while the continuous part corresponds to lower 
dimensional surfaces in the reference configuration where 
fractures or other type of nonzero dimensional singularities might take place.

\section{The radial problem}\label{rad:case}
For ease of exposition we limit ourselves in this section to the case where 
$m=3$. We recall that if $\tilde{W}$ is frame indifferent and isotropic then 
there is a symmetric function $\tilde{\Phi}$ such that
\begin{equation}\label{isotW}
\tilde{W}(\mathbf{F})=\tilde{\Phi} (v_1,v_2,v_3),
\end{equation}
where $v_1,v_2,v_3$ are the singular values of the matrix $\mathbf{F}$.
For the function $h(\cdot)$ in \eqref{h} we assume that it is strictly convex
so that it has a unique minimum at $d_0$, and that
\begin{equation}\label{growth:h:inf}
h(d)\sim Cd^\gamma,\quad d\To\infty,
\end{equation}
where $\gamma>1$ and $C$ is some positive constant.

For $\Omega$ equal to the unit ball with center at the origin, the radial
deformation
\begin{equation}\label{radmap}
\ts{u}(\ts{x})=\dfrac{r(R)}{R}\,\ts{x},\quad R=\norm{\ts{x}},
\end{equation}
has energy (up to a constant) given by:
\begin{equation}\label{rad_energy}
E_{\rm rad}(r)=\int_0^1R^2\left[\tilde{\Phi}\left(\deru{r}(R),\dfrac{r(R)}{R}
,\dfrac{r(R)}{R}\right)+h\left(\deru{r}(R)\left(\dfrac{r(R)}{R}
\right)^2\right)\right]\,\dif R.
\end{equation}
It is well known that for $p\in(1,3)$ in \eqref{growth}, there exists
$\lambda_c>d_0^{\frac{1}{3}}$ such that for $\lambda>\lambda_c$, the minimizer
$r_c$ of $E_{\rm rad}(\cdot)$ over the set
\begin{equation}\label{rad_adm}
\mc{A}_{\rm rad}=\set{r\in W^{1,1}(0,1)\,:\,\deru{r}(R)>0\mbox{ a.e.},
\,r(0)\ge0,\,r(1)=\lambda},
\end{equation}
exists and has $r_c(0)>0$.

With $v$ a radial function now, the modified functional \eqref{modfunct}
reduces up to a constant to:
\begin{eqnarray}
I_\eps^\tau(r,v)&=&\int_0^1R^2\left[\tilde{\Phi}\left(\deru{r}(R),\dfrac{r(R)}{R
}
,\dfrac{r(R)}{R}\right)+h\left(\deru{r}(R)\left(\dfrac{r(R)}{R}
\right)^2-v(R)\right)\right]\,\dif R\nonumber\\&&
+\int_0^1R^2\left[\frac{\eps^\alpha}{\alpha}|\deru{v}(R)|^\alpha
+\frac{1}{q\eps^q}\,\phi_\tau(v(R))\right]\,\dif R,\label{modfunctR}
\end{eqnarray}
and the set $\mc{U}$ becomes
\begin{eqnarray}
\mathcal{U}_{\rm rad} =\{(r,v)\in W^{1,1}(0,1)\times W^{1,\alpha}(0,1)
:r(0)=0,\ r(1)=\lambda,\nonumber\\  \deru{r}(R)(r(R)/R)^2>v(R)\ge0
\mbox{ a.e.},\ v(1)=0\}.\label{set:regR}
\end{eqnarray}
As a special case of Proposition \ref{affine_global_gen}, we now have the 
following result:
\begin{proposition}\label{affine_global}
Assume that the function $\tilde{\Phi}$ is quasi--convex. Then for $\lambda\le
d_0^{\frac{1}{3}}$, the global minimizer of $I_\eps^\tau(\cdot,\cdot)$ over
$\mathcal{U}_{\rm rad}$ is given by $r(R)=\lambda R$ and $v(R)=0$ for all $R$.
\end{proposition}
Note that if $(r,0)\in\mathcal{U}_{\rm rad}$, then $r(0)=0$, and
quasi--convexity implies that
\begin{equation}\label{qc:vzero}
I_\eps^\tau(r,0)\ge I_\eps^\tau(\lambda R,0).
\end{equation}
Moreover, since $I_\eps^\tau(r,0)=E_{\rm rad}(r)$, we have that
\begin{equation}\label{cav:vzero}
I_\eps^\tau(\lambda R,0)>E_{\rm rad}(r_c),\quad\lambda>\lambda_c.
\end{equation}
Using this we can now prove the following result.
\begin{theorem}\label{DPconvGrad}
Let $\gamma>1$ be as in \eqref{growth:h:inf} and assume that
the set of averages 
\begin{equation}\label{bddavg}
\set{\frac{1}{\tau}\int_0^\tau\phi_\tau(u)\,\dif 
u\,:\,\tau>\tau_0},
\end{equation}
is bounded for some $\tau_0>0$. For any $\lambda>\lambda_c$, there exists 
an $\eps_0>0$ such that for any $\eps\in(0,\eps_0)$, there exist $\tau(\eps)>0$ 
and a pair
$(\tilde{r}_\eps,\tilde{v}_\eps)
\in\mc{U}_{\rm rad}$ with $\tilde{v}_\eps$ non--constant, such that
\[
\inf_{\mc{U}_{\rm rad}}I_\eps^{\tau(\eps)}(r,v)\le 
I_\eps^{\tau(\eps)}(\tilde{r}_\eps,\tilde{v}_\eps)
<I_\eps^{\tau(\eps)}(\lambda R,0).
\]
Moreover $\tau(\eps)\To\infty$ as $\eps\To0^+$, and
\[
 \limite{\eps}{0}I_\eps^{\tau(\eps)}(\tilde{r}_\eps,\tilde{v}_\eps)=E_{\rm 
rad}(r_c).
\]
In particular, any minimizer $(r_\eps,v_\eps)$ must have $v_\eps$
non--constant, and
\begin{equation}\label{ineqRad}
\varliminf_{\eps\To0}\, I_\eps^{\tau(\eps)}(r_\eps,v_\eps)\le 
E_{\rm
rad}(r_c).
\end{equation}
\end{theorem}

\begin{prueba}
We now construct
$(\tilde{r},\tilde{v})$, $\tilde{v}$ non constant such that
\[
I_\eps^{\tau(\eps)}(\tilde{r},\tilde{v})<I_\eps^{\tau(\eps)}(\lambda R,0),
\]
for $\eps$ sufficiently small and $\tau(\eps)$ sufficiently large. For any 
$\delta>0$ 
we let
\begin{equation}\label{rest1}
\tau=\left(\frac{r_c(\delta)}{\delta}\right)^3-d_0.
\end{equation}
Since $r_c(0)>0$, we have that $\tau\To\infty$ as $\delta\To0^+$. For $\delta$ 
sufficiently small, we let $\eta\in(0,\delta)$ and define:
\[
\tilde{r}(R)=\left\{\begin{array}{rcl}\left[\frac{r_c(\delta)}{\delta} \right]
R&,&0\le R\le\delta,\\r_c(R)&,&\delta\le R\le1,\end{array}\right.
\]
\[
\tilde{v}(R)=\left\{\begin{array}{rcl}
\tau&,&0\le
R\le\delta-\eta,\\\frac{\tau}{\eta}(\delta-R)&, &\delta-\eta\le
R\le\delta,\\
0&,&\delta\le R\le1.\end{array}\right.
\]
For this test pair we have that
\begin{eqnarray*}
I_\eps^\tau(\tilde{r},\tilde{v})&=&
\int_0^{\delta-\eta}R^2\left[\tilde{\Phi}\left(\deru{\tilde{r}}(R),
\frac{\tilde{r}(R)}{R},\frac{\tilde{r}(R)}{R}\right)
+h\left(\deru{\tilde{r}}(R)\left[\frac{\tilde{r}(R)}{R}\right]^2
-\tilde{v}(R)\right)\right] \dif R\\
&&+ \int_{\delta-\eta}^{\delta}R^2\left[\tilde{\Phi}\left(\deru{\tilde{r}}(R),
\frac{\tilde{r}(R)}{R},\frac{\tilde{r}(R)}{R}\right)
+h\left(\deru{\tilde{r}}(R)\left[\frac{\tilde{r}(R)}{R}\right]^2
-\tilde{v}(R)\right)\right] \dif R\\
&&+ \int_{\delta}^{1}R^2\left[\tilde{\Phi}\left(\deru{\tilde{r}}(R),
\frac{\tilde{r}(R)}{R},\frac{\tilde{r}(R)}{R}\right)
+h\left(\deru{\tilde{r}}(R)\left[\frac{\tilde{r}(R)}{R}\right]^2
-\tilde{v}(R)\right)\right] \dif R\\
&&+\int_{\delta-\eta}^{\delta}R^2
\left[\frac{\eps^\alpha}{\alpha}\,|\deru{\tilde{v}}(R)|^\alpha+
\frac{1}{q\eps^q}\,\phi_\tau(\tilde{v}(R))\right]\dif R\equiv I_1+I_2+I_3+I_4.
\end{eqnarray*}
From the definition of $(\tilde{r},\tilde{v})$, it follows that:
\begin{enumerate}[i)]
\item
\begin{eqnarray*}
I_1&=&\int_0^{\delta-\eta}R^2\left[\tilde{\Phi}
\left(\frac{r_c(\delta)}{\delta},
\frac{r_c(\delta)}{\delta},\frac{r_c(\delta)}{\delta}\right)
+h(d_0)\right] \dif R\\
&=&\frac{(\delta-\eta)^3}{3} \left[\tilde{\Phi}
\left(\frac{r_c(\delta)}{\delta},
\frac{r_c(\delta)}{\delta},\frac{r_c(\delta)}{\delta}\right)
+h(d_0)\right].
\end{eqnarray*}
By taking
\begin{equation}\label{rest2}
\eta=\delta^{\beta_1},\quad\beta_1>1,
\end{equation}
we get from  that $I_1$ can be made arbitrarily small with 
$\delta$.
\item
For the term $I_2$, first note that since
\[
\tilde{v}(R)\le\tau=\left(\frac{r_c(\delta)}{\delta}\right)^3-d_0.
\]
we have that
\[
d_0\le\left(\frac{r_c(\delta)}{\delta}\right)^3-\tilde{v}(R)\le
\left(\frac{r_c(\delta)}{\delta}\right)^3.
\]
Since $h(\cdot)$ is increasing on $(d_0,\infty)$, it follows that
\[
h\left(\left(\frac{r_c(\delta)}{\delta}\right)^3-\tilde{v}(R)
\right)\le h\left(\left(\frac{r_c(\delta)}{\delta}\right)^3 \right).
\]
Thus
\begin{eqnarray*}
I_2&=&\int_{\delta-\eta}^{\delta}R^2\left[\tilde{\Phi}
\left(\frac{r_c(\delta)}{\delta},
\frac{r_c(\delta)}{\delta},\frac{r_c(\delta)}{\delta}\right)
+h\left(\left(\frac{r_c(\delta)}{\delta}\right)^3-\tilde{v}(R)
\right)\right] \dif R\\
&\le& \int_{\delta-\eta}^{\delta}R^2\left[\tilde{\Phi}
\left(\frac{r_c(\delta)}{\delta},
\frac{r_c(\delta)}{\delta},\frac{r_c(\delta)}{\delta}\right)
+h\left(\left(\frac{r_c(\delta)}{\delta}\right)^3
\right)\right] \dif R.
\end{eqnarray*}
Now
\[
\int_{\delta-\eta}^{\delta}R^2\tilde{\Phi}
\left(\frac{r_c(\delta)}{\delta},
\frac{r_c(\delta)}{\delta},\frac{r_c(\delta)}{\delta}\right)
\dif R\le\eta\delta^2 \tilde{\Phi}
\left(\frac{r_c(\delta)}{\delta},
\frac{r_c(\delta)}{\delta},\frac{r_c(\delta)}{\delta}\right).
\]
It follows from \eqref{growth} and \eqref{rest2} that the right hand side of the 
above inequality goes to zero with $\delta$. For the other term in $I_2$ we 
have:
\[
\int_{\delta-\eta}^{\delta}R^2h\left(\left(\frac{r_c(\delta)}{\delta}\right)^3
\right) \dif R\le C\frac{\eta\delta^2}{\delta^{3\gamma}},
\]
for some constant $C>0$ and where $\gamma>1$ is the growth rate of $h(d)$
as $d\To\infty$ (cf. \eqref{growth:h:inf}). If we further assume that 
$\beta_1>3\gamma-2$, then $I_2$ goes 
to zero with $\delta$.
\item
Since $\tilde{r}(R)=r_c(R)$ and $\tilde{v}(R)=0$ for $\delta\le R\le1$, we
have that
\begin{eqnarray*}
I_3&=& \int_{\delta}^{1}R^2\left[\tilde{\Phi}\left(\deru{r_c}(R),
\frac{r_c(R)}{R},\frac{r_c(R)}{R}\right)
+h\left(\deru{r_c}(R)\left[\frac{r_c(R)}{R}\right]^2
\right)\right] \dif R\\
&=&E_{\rm rad}(r_c)-
\int_{0}^{\delta}R^2\left[\tilde{\Phi}\left(\deru{r_c}(R),
\frac{r_c(R)}{R},\frac{r_c(R)}{R}\right)
+h\left(\deru{r_c}(R)\left[\frac{r_c(R)}{R}\right]^2
\right)\right] \dif R.
\end{eqnarray*}
But $R^2\left[\tilde{\Phi}\left(\deru{r_c}(R),
\frac{r_c(R)}{R},\frac{r_c(R)}{R}\right)
+h\left(\deru{r_c}(R)\left[\frac{r_c(R)}{R}\right]^2
\right)\right]\in L^1(0,1)$. Hence
\[
\int_{0}^{\delta}R^2\left[\tilde{\Phi}\left(\deru{r_c}(R),
\frac{r_c(R)}{R},\frac{r_c(R)}{R}\right)
+h\left(\deru{r_c}(R)\left[\frac{r_c(R)}{R}\right]^2
\right)\right] \dif R
\]
can be made arbitrarily small with $\delta$.
\item
For the last term in $I_\eps^\tau(\tilde{r},\tilde{v})$:
\begin{eqnarray*}
I_4&=& \int_{\delta-\eta}^{\delta}R^2
\left[\frac{\eps^\alpha}{\alpha}\,|\deru{\tilde{v}}(R)|^\alpha+
\frac{1}{q\eps^q}\,\phi_\tau(\tilde{v}(R))\right]\dif R\\
&\le&\eta\delta^2\left[C_1\frac{\eps^\alpha}{\delta^{3\alpha}\eta^\alpha}+\frac{
C_2 \eta} {\eps^q }
\right].
\end{eqnarray*}
Here we used that the set \eqref{bddavg} is bounded, and that
\[
 \int_{\delta-\eta}^\delta\phi_\tau(\tilde{v}(R))\,\dif
R=\frac{\eta}{\tau}\int_0^\tau\phi_\tau(u)\,\dif u\le C_2\eta.
\]
We set
\[
\frac{\eps^\alpha}{\delta^{3\alpha}\eta^{\alpha}}=\frac{\eta}{\eps^q} ,
\]
so that both terms on the right hand side of the inequality for $I_4$ above are 
of the same order, which upon recalling \eqref{rest2},
leads to
\begin{equation}\label{rest4}
\eps^{\alpha+q}=\delta^{(\beta_1+3)\alpha+\beta_1}.
\end{equation}
Given $\eps>0$, if $\delta$ is chosen according to \eqref{rest4}, then 
$\delta\To0^+$ and $\tau(\eps)\To\infty$ (cf. \eqref{rest1}) as $\eps\To0^+$. 
Thus
\[
\eta\delta^2\frac{\eta}{\eps^q}=\delta^{2\beta_1+2-\frac{q}{\alpha+q}
((\beta_1+3)\alpha+\beta_1)}=\delta^{\frac{\alpha}{\alpha+q}
(\beta_1-q)} ,
\]
and both terms in $I_4$ go to zero with $\delta$ (and thus with $\eps$) provided
$\beta_1>q$.
\end{enumerate}
Thus we can conclude that 
\[
I_\eps^{\tau(\eps)}(\tilde{r},\tilde{v})\To E_{\rm rad}(r_c),\quad\mbox{as 
}\eps\To0,
\]
and from \eqref{cav:vzero}, that $I_\eps^{\tau(\eps)}(\lambda 
R,0)>I_\eps^{\tau(\eps)}(\tilde{r},\tilde{v})$
for $\eps$ sufficiently small.
\end{prueba}

Now, in the radial case, the limiting function $\ts{u}^*$ of Theorem 
\ref{reppnumsch} must be radial, and the limiting measure $\nu$ must be a 
non--negative multiple of the Dirac delta distribution centered at the origin.
Since $\ts{u}^*$ is radial we must have, with $\Omega$ the unit ball, that 
\[
 \int_\Omega W(\nabla\ts{u}^*)\,\dif\ts{x}\ge E_{\rm rad}(r_c).
\]
Thus, combining this with \eqref{ineqG2} and \eqref{ineqRad}, we get
that $P(B,\Omega)=0$, and that 
\[
E_{\rm rad}(r_c)=\int_\Omega 
W(\nabla\ts{u}^*)\,\dif\ts{x}=\varliminf_{\eps\To0}\, \inf_{\mc{U}_{\rm 
rad}}I_\eps^{\tau(\eps)}(r,v),
\]
where $\ts{u}^*$ is given by \eqref{radmap} using $r_c$. Thus we have proved the 
following:
\begin{theorem}\label{convgRad}
Fix $\lambda>\lambda_c$ and let $(r_\eps^\tau,v_\eps^\tau)$ be a minimizer of 
$I_\eps^\tau$ 
over $\mathcal{U}_{\rm rad} $ and let $\ts{u}_\eps^\tau$ be the radial 
map \eqref{radmap} corresponding to $r_\eps^\tau$. Let $\set{\eps_j}$ be a 
sequence such that $\eps_j\To0$. Then for a subsequence of $\set{\eps_j}$, 
there exists a sequence $\set{\tau_j}$ with $\tau_j\To\infty$, such that the 
sequences $\set{\ts{u}_j}$ and $\set{v_j}$, where 
$\ts{u}_j=\ts{u}_{\eps_j}^{\tau_j}$ and $v_j=v_{\eps_j}^{\tau_j}$, have
subsequences $\set{\ts{u}_{j_k}}$ and $\set{v_{j_k}}$ with
$\ts{u}_{j_k}\rightharpoonup\ts{u}^*$ in $W^{1,p}(\Omega)$ and
$\wsconv{v_{j_k}}{\nu}$ in $\mc{M}(\Omega)$, where $\ts{u}^*$ is given by 
\eqref{radmap} using $r_c$ (the minimizer of $E_{\rm rad}(\cdot)$ over the set
\eqref{rad_adm}) and $\nu=\kappa\delta_{\ts{0}}$ with $\kappa>0$. Moreover
\[
E_{\rm 
rad}(r_c)=\varliminf_{k\To\infty}I_{\eps_{j_k}}^{\tau_{j_k}}(r_{j_k},v_{j_k}).
\]
\end{theorem}

\subsection{The Euler--Lagrange equations}
In this section we show that the minimizers of \eqref{modfunctR} over 
\eqref{set:regR}, satisfy the Euler--Lagrange equations for this 
functional. The analysis is not straightforward, basically due to the singular 
behaviour of the function $h(\cdot)$ (cf. \eqref{h}), and the inequality 
constraints involving the phase function $v$, that is, its non--negativity and 
the inequality involving the determinant of the deformation $r$. The proof is a 
variation of that in \cite{Ba82}.

For the following discussion we use the notation:
\begin{equation}\label{phi4}
 \hat{\Phi}(v_1,v_2,v_3,v_4)=\tilde{\Phi}(v_1,v_2,v_3)+h(v_1v_2v_3-v_4).
\end{equation}
Also we shall write
\[
 \hat{\Phi}(r(R),v(R))=\hat{\Phi}\left(\deru{r}(R),\frac{r(R)}{R},\frac{r(R)}{R
} ,v(R)\right),\quad\mbox{etc.}
\]
The functional \eqref{modfunctR} can now be written as:
\begin{eqnarray}
I_\eps^\tau(r,v)&=&\int_0^1R^2\hat{\Phi}\left(r(R),v(R)\right) \dif 
R\nonumber\\~~~~~~~~~~~~~~~~~~&&+
\int_0^1R^2\left[\frac{\eps^\alpha}{\alpha}|\deru{v}(R)|^\alpha
+\frac{1}{q\eps^q}\phi_\tau(v(R))\right]\,\dif R,\label{modfunctR2}
\end{eqnarray}
where $(r,v)\in\mathcal{U}_{\rm rad}$ (cf. \eqref{set:regR}). 

For the analysis in this section we take $\tilde{\Phi}$ in \eqref{phi4} as
 \begin{equation}\label{Ogdentype}
  \tilde{\Phi}(v_1,v_2,v_3)=\sum_{i=1}^{3}\psi(v_i),
 \end{equation}
 where $\psi$ is a non-negative convex $C^3$ function over $(0,\infty)$, 
and with 
\begin{equation}\label{const_H1}
 \abs{v\,\deru{\psi}(cv)}\le K\psi(v),\quad\abs{c-1}\le\gamma_0,
\end{equation}
for all $v>0$ and for some positive constants $K,\gamma_0$. 
However, our results hold as well for more general stored energy functions 
under suitable assumptions. We now have:
\begin{theorem}\label{mod_EL1}
 Let $(r,v)$ be any minimizer of $I_\eps^\tau$ over \eqref{set:regR}. 
Assume that the
functions $h(\cdot)$ and $\psi(\cdot)$ in \eqref{phi4} together with 
\eqref{Ogdentype}, satisfy \eqref{h} and \eqref{const_H1} respectively.
Then $(r,v)\in
C^1(0,1]\times C^1(0,1]$, $\deru{r}(R)>0$ for all $R\in(0,1]$, 
$R^{2}\hat{\Phi}_1(r(R),v(R))$ is
$C^1(0,1]$, and
\begin{subequations}\label{EL:interf}
\begin{eqnarray}
&\displaystyle\td{}{R}\left[R^2\hat{\Phi}_{,1}(r(R),v(R))\right]=
2R\hat{\Phi}_{,2}(r(R),v(R)),\quad0<R<1,&\label{ELeqr1}\\
&\displaystyle
v^{\frac{1}{2}}(R)\bigg(\eps^\alpha\td{}{R}[R^2|\deru{v}(R)|^{\alpha-1}
\sgn(\deru{v}(R))]
~~~~~~~~~~~~~~~~~~~~~~~~~~~~~~~~~~~~~~~~~~~~\nonumber&\\
&-R^2\left[\hat{\Phi}_{,4}(r(R),v(R))+\frac{1}{q\eps^q}\deru{\phi_\tau}
(v(R))\right]\bigg)
=0,\quad0<R<1,&\label{ELeqv}
\end{eqnarray}
\end{subequations}
with boundary conditions:
\begin{equation}\label{EL:interfBC}
r(0)=0,\quad r(1)=\lambda,\quad \limite{R}{0^+}
R^2|\deru{v}(R)|^{\alpha-1}\sgn(\deru{v}(R))v^{\frac{1}{2}}(R)=0,\quad v(1)=0.
\end{equation}
\end{theorem}
\begin{prueba}
If we let 
$v=u^2$, then our problem is equivalent to that of minimizing
\begin{eqnarray}
\hat{I}_\eps^\tau(r,u)&=&\int_0^1R^2\hat{\Phi}\left(r(R),
u^2(R)\right) \dif R\nonumber\\~~~~~~~~~~~~~~~~~~&&+
\int_0^1R^2\left[\frac{\eps^\alpha}{\alpha}|2u(R)\deru{u}(R)|^\alpha
+\frac{1}{q\eps^q}\phi_\tau(u^2(R))\right]\,\dif R,\label{modfunctR2u}
\end{eqnarray}
over
\begin{eqnarray}
\hat{\mathcal{U}}_{\rm rad} =\{(r,u)\in W^{1,1}(0,1)\times W^{1,\alpha}(0,1)
:r(0)=0,\ r(1)=\lambda,\nonumber\\  \deru{r}(R)(r(R)/R)^2>u^2(R)\mbox{ a.e.}
,\ u(1)=0\}.\label{set:regRu}
\end{eqnarray}
Note that since $u\in W^{1,\alpha}(0,1)$, then $u$ is continuous in $[0,1]$. 
Hence 
both $u^2$ and $u\deru{u}$ belong to $L^\alpha(0,1)$.

Let $(r,u)$ be any minimizer of $\hat{I}_\eps^\tau$ over \eqref{set:regRu}. 
We first consider variations only in $r$, keeping $u$ fixed. We make the change 
of variables $w=r^3(R)$ and $\rho=R^3$. It follows now that
\[
 \dot{w}(\rho)=\td{w}{\rho}(\rho)=\deru{r}(R)\left(\frac{r(R)}{R}\right)^2.
\]
The first part of the functional \eqref{modfunctR2u} can now be written as
\[
 \int_0^1f(\rho,w,\dot{w},u^2)\,\dif\rho,
\]
where
\[
3f(\rho,w,\dot{w},u^2)=\tilde{\Phi}((\rho/w)^{\frac{2}{3}}\dot{w}
,(w/\rho)^{\frac{1}{3}},(w/\rho)^{\frac{1}{3}})+h(\dot{w}-u^2).
\]
For 
$k\ge1$ we define
\[
S_k=\set{\rho\in\left(\frac{1}{k},1\right)\,:\,\frac{1}{k}\le
\dot{w}(\rho)-u^2(\rho)\le k},
\]
and let $\chi_k$ its characteristic function. Let $\omega\in 
L^\infty(0,1)$ be such that 
\[
\int_{S_k} \omega(s)\,\dif s=0, 
\]
and for any $\gamma>0$, define the variations
\[
 w_\gamma(\rho)=w(\rho)+\gamma\int_0^\rho \chi_k(s)\omega(s)\,\dif s.
\]
Note that $w_\gamma(0)=0$ and $w_\gamma(1)=\lambda^3$. The rest of the proof,  
using \eqref{const_H1}, is as in \cite{Ba82}, from which it follows (after 
changing back to $R$ and $r$) that $r\in C^1(0,1]$, $\deru{r}(R)>0$ for all 
$R\in(0,1]$, $R^{2}\hat{\Phi}_1(r(R),u(R))$ is $C^1(0,1]$, and that 
equations \eqref{ELeqr1} and the first two boundary conditions in 
\eqref{EL:interfBC} hold. 

We now consider variations in $u$ keeping $r$ fixed. For any $k\ge1$, let $z\in 
W^{1,\infty}(0,1)$ have support in $(\frac{1}{k},1)$. and let
\[
 u_\gamma=u+\gamma z.
\]
Note that $u_\gamma(1)=0$. Moreover, since $r\in C^1(0,1]$ and $u\in C[0,1]$, 
it 
follows that
\[
 \deru{r}(R)\left(\frac{r(R)}{R}\right)^2>u_\gamma^2(R),\quad 
R\in\left[\frac{1}{k},1\right],
\]
for $\gamma$ sufficiently small. It follows now, upon setting 
$\delta(R)=\deru{r}(R)(r(R)/R)^2$, that
\begin{eqnarray*}
\dfrac{\hat{I}_\eps^\tau(r,u_\gamma)-\hat{I}_\eps^\tau(r,u)}{\gamma}
&=&\frac{1}{\gamma}\int_0^1R^2\left[h(\delta-u_\gamma^2)
-h(\delta-u^2)\right]\dif R\\&&+\frac{1}{\gamma}
\int_0^1\frac{\eps^\alpha}{\alpha}R^2\left[|2u_\gamma 
\deru{u_\gamma}|^\alpha-|2u\deru{u}|^\alpha
\right]\,\dif R\\&&+\frac{1}{\gamma}
\int_0^1\frac{1}{q\eps^q}R^2\left[\phi_\tau(u_\gamma^2)
-\phi_\tau(u^2)\right]\,\dif R.
\end{eqnarray*}
Now
\begin{multline*}
\displaystyle \frac{1}{\gamma}\int_0^1R^2\left[h(\delta-u_\gamma^2)
-h(\delta-u^2)\right]\dif R=\\
\displaystyle
\frac{1}{\gamma}\int_0^1R^2\int_0^1\td{}{t}[
h(\delta-(tu_\gamma^2+(1-t)u^2))]\,\dif t\dif R=\\
\displaystyle
-\int_{\frac{1}{k}}^1R^2z(2u+\gamma z)\int_0^1
\deru{h}(\delta-(tu_\gamma^2+(1-t)u^2))\,\dif t\dif R\\
\To\displaystyle-\int_{\frac{1}{k}}^12
\deru{h}(\delta-u^2)uzR^2\dif R,
\end{multline*}
as $\gamma\To0$. Similarly
\begin{align*}
 \frac{1}{\gamma}
\int_0^1\frac{\eps^\alpha}{\alpha}R^2\left[|2u_\gamma 
\deru{u_\gamma}|^\alpha-|2u\deru{u}|^\alpha
\right]\,\dif R&\To \int_{\frac{1}{k}}^1\eps^\alpha 
|2u\deru{u}|^{\alpha-1}
\sgn(2u\deru{u})2\deru{(uz)}
R^2\,\dif R,\\
 \frac{1}{\gamma}
\int_0^1\frac{1}{q\eps^q}R^2\left[\phi_\tau(u_\gamma^2)
-\phi_\tau(u^2)\right]\,\dif R&\To
\int_{\frac{1}{k}}^1\frac{1}{q\eps^q}\deru{\phi_\tau}(u^2)
2uzR^2\,\dif R,
\end{align*}
as $\gamma\To0$. Since
\[
 \limite{\gamma}{0}
 \dfrac{\hat{I}_\eps^\tau(r,u_\gamma)-\hat{I}_\eps^\tau(r,u)}{\gamma}=0,
\]
we get, combining our previous results that
\[
 \int_{\frac{1}{k}}^1[-\deru{h}(\delta-u^2)uz+
 \eps^\alpha |2u\deru{u}|^{\alpha-1}
\sgn(2u\deru{u})\deru{(uz)}+
\frac{1}{q\eps^q}\deru{\phi_\tau}(u^2)
uz]R^2\dif R=0,
\]
or after collecting terms,
\begin{multline*}
  \int_{\frac{1}{k}}^1[\eps^\alpha |2u\deru{u}|^{\alpha-1}
\sgn(2u\deru{u})u\deru{z}
 +(\eps^\alpha |2u\deru{u}|^{\alpha-1}
\sgn(2u\deru{u})\deru{u}+\\
\frac{1}{q\eps^q}\deru{\phi_\tau}(u^2)
u-\deru{h}(\delta-u^2)u)z]R^2\dif R=0,
\end{multline*}
for all $z\in W^{1,\infty}(0,1)$ with support in $(\frac{1}{k},1)$. The 
coefficient of $z$ in this expression is in $L^1(\frac{1}{k},1)$. Hence the 
above equation is equivalent to
\begin{multline*}
  \int_{\frac{1}{k}}^1\left[\eps^\alpha |2u\deru{u}|^{\alpha-1}
\sgn(2u\deru{u})uR^2+\int_R^1
 (\eps^\alpha |2u\deru{u}|^{\alpha-1}
\sgn(2u\deru{u})\deru{u}+\right.\\\left.
\frac{1}{q\eps^q}\deru{\phi_\tau}(u^2)
u-\deru{h}(\delta-u^2)u)\xi^2\,\dif\xi\right]\deru{z}\dif R=0.
\end{multline*}
The arbitrariness of $z$ implies now that for some constant $C$ independent of 
$k$, we have
\begin{multline*}
  \eps^\alpha |2u\deru{u}|^{\alpha-1}
\sgn(2u\deru{u})uR^2+\int_R^1
 (\eps^\alpha |2u\deru{u}|^{\alpha-1}
\sgn(2u\deru{u})\deru{u}+\\
\frac{1}{q\eps^q}\deru{\phi_\tau}(u^2)
u-\deru{h}(\delta-u^2)u)\xi^2\,\dif\xi=C,
\end{multline*}
over $(0,1)$. It follows from this equation that over the intervals where 
$u\ne0$, the function $|2u\deru{u}|^{\alpha-1}
\sgn(2u\deru{u})R^2$ is absolutely continuous. Hence after 
differentiating and simplifying, the equation above yields that
\[
  \left(\eps^\alpha \td{}{R}\left[|2u\deru{u}|^{\alpha-1}
\sgn(2u\deru{u})R^2\right]
-\left(\frac{1}{q\eps^q}\deru{\phi_\tau}(u^2)
-\deru{h}(\delta-u^2)\right)R^2\right)u=0,
\]
i.e., that \eqref{ELeqv} holds after reverting the substitution $v=u^2$. A 
standard argument now using variations $z$ not vanishing at $R=0$, yields the 
third boundary condition in \eqref{EL:interfBC}.
\end{prueba}
\begin{remark}
Note that the pair $r(R)=\lambda R$ and $v(R)=0$ is a solution of 
\eqref{EL:interf}-\eqref{EL:interfBC} for all $\lambda$. By Proposition 
\ref{affine_global}, this pair is a global minimizer for 
$\lambda<d_0^{\frac{1}{3}}$. However for $\lambda>\lambda_c$, $\eps$ 
sufficiently small, and $\tau$ sufficiently large, we get from Theorems 
\ref{DPconvGrad} and \ref{mod_EL1}, that the minimizer must have $v$ 
non--constant, with segments in which $v$ vanishes, and (non--trivial) 
segments in which the differential equation
\[
\eps^\alpha\td{}{R}[R^2|\deru{v}(R)|^{\alpha-1}
\sgn(\deru{v}(R))]
=R^2\left[\hat{\Phi}_{,4}(r(R),v(R))+\frac{1}{q\eps^q}\deru{\phi_\tau}
(v(R))\right],
\]
holds.
\end{remark}

\subsection{Numerical results}

To approximate the minimum of \eqref{modfunctR} over \eqref{set:regR}, let 
$h=1/m$ and $R_i=ih$, $0\le i\le m$, where $m\ge1$. We write $(r_i,v_i)$ for 
any approximation of $(r(R_i,v(R_i)))$, $0\le i\le m$, and 
\[
 R_{i-\frac{1}{2}}=\frac{R_i+R_{i-1}}{2},\quad \delta 
r_{i-\frac{1}{2}}=\frac{r_i-r_{i-1}}{h},\quad 
\left(\frac{r}{R}\right)_{i-\frac{1}{2}}=\frac{r_i+r_{i-1}}{R_i+R_{i-1}},
\quad  i=1,\ldots,m.
\]
Now we discretize $I_\eps^\tau$ as follows:
\begin{multline}
I_{\eps,h}^\tau=
h\sum_{i=1}^{m}R_{i-\frac{1}{2}}^{2}\left[\tilde{\Phi}\left(\delta 
r_{i-\frac{1}{2}},\left(\frac{r}{R}\right)_{i-\frac{1}{2}}, 
\left(\frac{r}{R}\right)_{i-\frac{1}{2}}\right)+h\left(
\delta r_{i-\frac{1}{2}}\left(\frac{r}{R}\right)_{i-\frac{1}{2}}^2
-v_{i-\frac{1}{2}}\right) \right]\\
+h\sum_{i=1}^{m}R_{i-\frac{1}{2}}^{2}\left[\frac{\eps^\alpha}{\alpha} 
\abs{\delta 
v_{i-\frac{1}{2}}}^\alpha+\frac{1}{q\eps^q}\,\phi_\tau(v_{i-\frac{1}{2}}) 
\right],\label{discRfunc}
\end{multline}
subject to $r_0=0$, $r_m=\lambda$, $v_m=0$ and 
\begin{equation}\label{discRconst}
 v_i\ge0,\quad0\le i\le m,\quad \delta 
r_{i-\frac{1}{2}}\left(\frac{r}{R}\right)_{i-\frac{1}{2}}^2
-v_{i-\frac{1}{2}}>0,\quad1\le i\le m.
\end{equation}
We compute (relative) minimizers of \eqref{discRfunc} over 
\eqref{discRconst} using the function \texttt{fmincon} of MATLAB 
 with the option for an interior point algorithm. With this routine the first 
set of conditions in \eqref{discRconst} can be directly specified as lower 
bounds on 
the $v_i$'s, while the second set of constraints is specified with the 
option for inequality constraints. The strict sign in the second set of 
conditions in \eqref{discRconst} is indirectly handled by the interior 
point algorithm with the $h$ playing the role of an interior penalty function 
(since $h(d)\To\infty$ as $d\searrow0$). For the various functions in the 
functional above we used the following:
\begin{eqnarray}
\tilde{\Phi}(v_1,v_2,v_3)&=&\mu\left(v_1^p+v_2^p+v_3^p\right),\quad
h(d)= c_1d^\gamma+c_2d^{-\delta},\label{numPhih}\\
\phi_\tau(v)&=&\left\{\begin{array}{ccl}Kv^2(v-\tau)^2&,&v\in[0,\tau],\\ 
0&,&\mbox{elsewhere,}\end{array}\right.\label{numphi}
\end{eqnarray}
where $p\in[1,3)$, $\mu,c_1,c_2\ge0$, $\gamma,\delta\ge1$, and $K>0$ is chosen 
so that $\max_{v\in[0,\tau]}\phi_\tau(v)=1$. With this choice of $K$ the 
boundedness condition \eqref{bddavg} holds as well. In the 
calculations below we use 
$m=100$ and the following values for the various constants:
\begin{equation}\label{numdata}
\mu=1.0,\quad c_1=1.0,\quad p=2.0, \quad\alpha=2,\quad
\gamma=2.0,\quad\delta=2.0,\quad \tau=3,\quad\lambda=1.5,
\end{equation}
with $c_2=(p\mu+\gamma c_1)/\delta$ so as to make the reference configuration 
stress free. In this case the minimizer $r_c$ of \eqref{rad_energy} over 
\eqref{rad_adm} has $E_{\rm rad}(r_c)\approx4.5396$ with $r_c(0)\approx1.222$, 
while the affine deformation $r^h(R)=\lambda R$ has energy $E_{\rm 
rad}(r^h)\approx 6.1053$.
\begin{table}[t]
\begin{center}
\begin{tabular}{|c|c|c|}\hline
$\eps^2$ & $I_{\eps,h}^\tau$&$\delta r_{\frac{1}{2}}$\\\hline\hline
$10^{-5}$& $6.101645$&$5.412$\\\hline
$10^{-6}$& $6.105267$&$1.560$\\\hline
$10^{-7}$& $6.105291$&$1.590$\\\hline
$10^{-8}$& $5.634048$&$32.85$\\\hline
$10^{-9}$& $4.771748$&$50.47$\\\hline
$10^{-10}$& $4.535530$&$49.91$\\\hline
\end{tabular}
\end{center}
\caption{Convergence of the decoupled penalized scheme in the 
radial case using \eqref{numPhih} and \eqref{numphi} with data 
\eqref{numdata}.}\label{tab:1}
\end{table}
\begin{figure}[ht]
\begin{center}
\scalebox{0.6}{\includegraphics{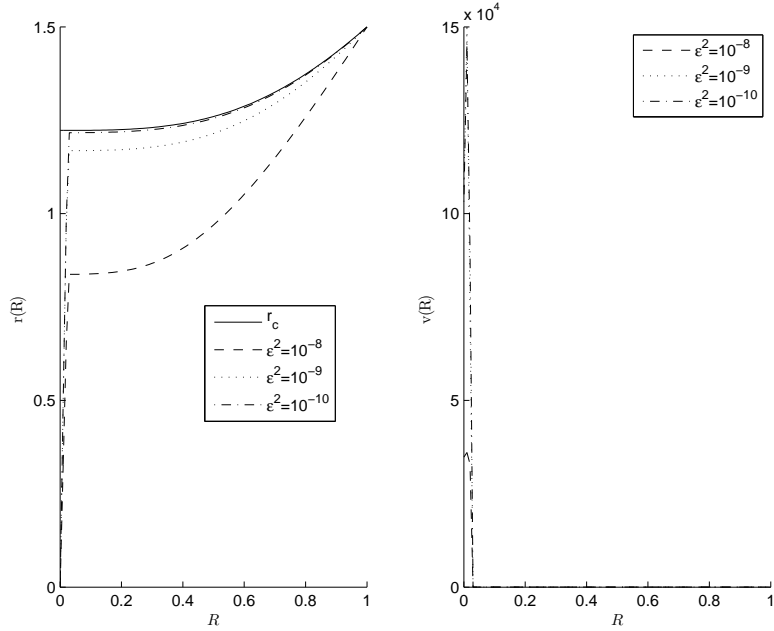}}
\end{center}
\caption{Numerical results for the data in \eqref{numdata}.}
\label{fig:1.2}
\end{figure}

In Table \ref{tab:1} we show the computed minimum energies for different values 
of $\eps^2$. In each case the 
iterations were started from the discretized 
versions of the affine deformation $r^h$ and $v=0$. From the values in the 
table we see that the approximations of $r$ for
$\eps^2=10^{-5},10^{-6},10^{-7}$ stay ``close'' to the affine deformation $r^h$ 
but developing a steep slope close to $R=0$. This process picks up after 
$\eps^2=10^{-8}$, where the energies get very close to the energy $E_{\rm 
rad}(r_c)\approx4.5396$ of the cavitated solution, and with very large slopes 
close to $R=0$. In Figure \ref{fig:1.2} (left) we show the computed $r$ 
approximations for $\eps^2=10^{-8},10^{-9},10^{-10}$ which are clearly 
converging to the cavitated solution $r_c$. On the other hand, Figure 
\ref{fig:1.2} (right) shows the corresponding approximations of $v$, 
which are clearly developing a singularity close to $R=0$ to match the 
corresponding singular behaviour of the determinants of the $r$ approximations.

\section{An elastic fluid}\label{elast:case}
We now consider the case of an \textit{elastic fluid} in which $\tilde{W}=0$ in 
\eqref{modfunct}. We will
show that in this case, if $\Omega$ is the unit ball with centre at the origin,
then any minimizer $(\ts{u},v)\in\mc{U}$ of $I$ must be radial
with the phase function $v$ radially decreasing. The main concept to get this
result is that of \textit{function rearrangement}.

For any set $A\subset\Real^3$, its \textit{symmetric rearrangement} $A^*$
is given by the ball with center at the origin and with measure equal to
that of $A$. Given $g:\Omega\To\Real$ a non--negative measurable such that
\[
\mbox{meas}\set{\ts{x}\in\Omega:g(\ts{x})>t}<\infty,\quad\forall\,t>0,
\]
we define the \textit{symmetric decreasing rearrangement} $g^*$ of $g$ by:
\[
g^*(\ts{x})=\int_0^\infty\chi_{\set{\ts{s}\,:\,g(\ts{s})>t}^*}(\ts{x})\,\dif t.
\]
The symmetric decreasing rearrangement of a function have several important
properties (see eg. \cite{LiLo2001}) that we list here without proof:
\begin{enumerate}
\item[R1:]
$g^*$ is radially symmetric and decreasing, and 
$\mbox{supp}(g^*)\subset\Omega^*$.
\item[R2:]
$\mbox{meas}\set{\ts{x}\in\Omega:g(\ts{x})>t}=
\mbox{meas}\set{\ts{x}\in\Omega^*:g^*(\ts{x})>t}$ for
all $t\ge0$.
\item[R3:]
For any function $f$ that is the difference of two monotone functions,
\[
\int_{\Omega}f(g(\ts{x}))\,\dif 
\ts{x}=\int_{\Omega^*}f(g^*(\ts{x}))\,\dif\ts{x}.
\]
In particular,
\[
\int_{\Omega}g(\ts{x})\,\dif \ts{x}=\int_{\Omega^*}g^*(\ts{x})\,\dif\ts{x}.
\]
\item[R4:]
If $g\in W^{1,p}(\Omega)$, $p\ge1$, then
\[
\int_{\Omega^*}||\nabla g^*(\ts{x})||^p\,\dif\ts{x}\le
\int_{\Omega}||\nabla g(\ts{x})||^p\,\dif\ts{x}.
\]
\end{enumerate}

Let $\Omega=\mc{B}$ be the unit ball with center at the origin. With 
$\tilde{W}=0$, our functional \eqref{modfunct} reduces to:
\begin{equation}\label{modfunctPF}
I(\ts{u},v)=\int_{\mc{B}}
h(\det\nabla\ts{u}
(\ts{x})-v(\ts{x}))\,\dif\ts{x}
+\int_{\mc{B}}\left[\frac{\eps^\alpha}{\alpha} ||\nabla
v(\ts{x})||^\alpha+\frac{1}{q\eps^q}\phi_\tau(v(\ts{x}))\right]
\,\dif\ts{x}.
\end{equation}
For $(\ts{u},v)\in\mc{U}$, let
\[
w(\ts{x})=\det\nabla\ts{u}(\ts{x})-v(\ts{x}),\quad\ts{x}\in\mc{B}.
\]
Note that
\[
\int_{\mc{B}}w(\ts{x})\,\dif\ts{x}=\int_{\mc{B}}\det\nabla\ts{u}(\ts{x})\,
\dif\ts{x}-
\int_{\mc{B}}v(\ts{x})\,\dif\ts{x}=\frac{4\pi}{3}\lambda^3- 
\int_{\mc{B}}v(\ts{x})\,\dif\ts{x}.
\]
Thus instead of minimizing \eqref{modfunctPF} over $\mc{U}$, we consider the 
problem of minimizing
\begin{equation}\label{modfunctPF2}
J(w,v)=\int_{\mc{B}}
h(w(\ts{x}))\,\dif\ts{x}
+\int_{\mc{B}}\left[\frac{\eps^\alpha}{\alpha} ||\nabla
v(\ts{x})||^\alpha+\frac{1}{q\eps^q}\phi_\tau(v(\ts{x}))\right]
\,\dif\ts{x}.
\end{equation}
over the set
\begin{eqnarray} \label{set:reg2}
\mathcal{J} =\{(w,v)\in L^{p}(\Omega)\times W^{1,\alpha}(\Omega)
:\ \int_{\mc{B}}(w+v)(\ts{x})\,\dif\ts{x}=\frac{4\pi}{3}\lambda^3,\\ \nonumber
w>0\ {\rm a.e.},\ v\ge0\ {\rm a.e.},\ \ v|_{\partial\mc{B}}=0\}
\end{eqnarray}
We now have:
\begin{proposition}\label{PF1}
Let
\begin{equation}\label{hmonotD}
h(d)=c_1d^\gamma+c_2d^{-\delta},
\end{equation}
where $c_1$ and $c_2$ are positive, $\gamma>1$, $\delta>0$. Then any minimizer 
$(w,v)$ of
$J(\cdot,\cdot)$ over $\mc{J}$ must have $w,v$ both radially symmetric with $v$
radially decreasing.
\end{proposition}

\begin{prueba}
For any $(w,v)\in\mc{J}$, we consider its rearrangement $(w^*,v^*)$. Clearly 
this pair has the properties of been radially symmetric with $v^*$ radially 
decreasing. We now check that $(w^*,v^*)\in\mc{J}$:
\begin{enumerate}[i)]
\item
By construction $w^*,v^*\ge0$ and since (property R2 with $t=0$)
\[
\mbox{meas}\set{\ts{x}\in\mc{B}:w^*(\ts{x})>0}= 
\mbox{meas}\set{\ts{x}\in\mc{B}:w(\ts{x})>0},
\]
it follows that $w^*>0$ a.e.
\item
Since $\mbox{supp}(v^*)\subset\mc{B}^*=\mc{B}$, it follows that 
$v^*|_{\partial\mc{B}}=0$.
\item
To check the integral constraint, we note that by the property R3 of
rearrangements listed above,
\[
\int_{\mc{B}}(w^*+v^*)(\ts{x})\,\dif\ts{x}= 
\int_{\mc{B}}(w+v)(\ts{x})\,\dif\ts{x}=\frac{4\pi}{3}\lambda^3.
\]
\end{enumerate}
Thus $(w^*,v^*)\in\mc{J}$. Finally, from the properties R3 (with $f=h$) and R4
listed above of the rearrangements, we get that
\[
J(w^*,v^*)\le J(w,v),
\]
from which the result follows.
\end{prueba}
\begin{theorem}
Let \eqref{hmonotD} holds. Then any minimizer $(\ts{u},v)\in\mc{U}$ of
\eqref{modfunctPF} must be radially symmetric with the phase function $v$
radially decreasing.
\end{theorem}
\begin{prueba}
We show that for any $(\ts{u},v)\in\mc{U}$, there exists a pair
$(\hat{\ts{u}},\hat{v})\in\mc{U}$ radially symmetric with $\hat{v}$ radially
decreasing such that $I(\hat{\ts{u}},\hat{v})\le I(\ts{u},v)$.

For any $(\ts{u},v)\in\mc{U}$, the pair $(w,v)\in\mc{J}$ where 
$w=\det\nabla\ts{u}-v$. By Proposition \ref{PF1} we get that there exists a 
pair 
$(\hat{w},\hat{v})\in\mc{J}$ radially symmetric with $\hat{v}$ radially 
decreasing such that
\[
J(\hat{w},\hat{v})\le J(w,v)=I(\ts{u},v).
\]
Define
\[
r^3(R)=\int_0^Rs^2(\hat{w}+\hat{v})(s)\,\dif s,\quad R\in[0,1],
\]
and
\[
\hat{\ts{u}}(\ts{x})=r(||\ts{x}||)\,\frac{\ts{x}}{||\ts{x}||},\quad\ts{x}\in\mc{
B}.
\]
It follows now that 
$\det\nabla\hat{\ts{u}}(\ts{x})=(\hat{w}+\hat{v})(||\ts{x}||)$ and that
\[
r^3(1)=\int_0^1s^2(\hat{w}+\hat{v})(s)\,\dif 
s=\frac{3}{4\pi}\int_{\mc{B}}(\hat{w}+\hat{v})(||\ts{x}||)\,\dif\ts{x}
=\lambda^3.
\]
This together with the other properties of $(\hat{w},\hat{v})$ gives that 
$(\hat{\ts{u}},\hat{v})\in\mc{U}$. Moreover
\[
I(\hat{\ts{u}},\hat{v})=J(\hat{w},\hat{v})\le J(w,v)=I(\ts{u},v),
\]
from which the result follows.
\end{prueba}

It follows now that we can restrict the minimization in \eqref{modfunctPF} to 
radial functions. Thus from Theorem \ref{convgRad} we get that that the 
corresponding minimizers and minimizing energies given by \eqref{modfunctPF}, 
converge to the radial cavitating minimizer corresponding to $r_c$ and its 
corresponding energy.

\subsection{Numerical results}
We now describe some numerical simulations using a gradient flow 
iteration together with the finite element method, to minimize the functional 
\eqref{modfunctPF} directly. The inequality $\det\nabla\bu>v$ in $\mc{U}$ (cf. 
\eqref{set:reg}) is handled numerically via an interior point 
penalization technique, using the the singular behaviour in $h(\cdot)$ as the 
penalization function. The other inequality, namely, $v\ge0$, is handled by 
extending the function $\phi_\tau$ to negative arguments, with a large 
penalization for these values. In particular we use
\begin{equation}\label{extphi}
 \phi_\tau(v)=\left\{\begin{array}{ccl}Kv^2(v-\tau)^2&,&v\in[0,\tau],\\ 
0&,&v>\tau,\\Mv^2&,&v<0,\end{array}\right.
\end{equation}
where $M>0$ is a large parameter.

For simplicity we take $\alpha=q=2$ in \eqref{modfunctPF}. It follows 
now formally, that the Euler--Lagrange equations for \eqref{modfunctPF} are  
given by:
\begin{subequations}\label{pdeEF}
\begin{eqnarray} 
-\mbox{div}\left[\deru{h}(\det\nabla\bu-v)(\mbox{adj}\nabla\bu)^T\right]=\ts{0},
\label{pde1EF}\\
-\eps^2\Delta v-\deru{h}(\det\nabla\bu-v)
+\frac{1}{2\eps^2}\deru{\phi_\tau}(v)=0,\label{pde2EF}
\end{eqnarray}
\end{subequations}
in $\mc{B}$, with boundary conditions $\bu=\lambda \ts{I}$ and $v=0$ on 
$\partial\mc{B}$.

A gradient flow iteration (cf. \cite{Neu1997}) to minimize \eqref{modfunctPF}, 
assumes that $\bu$ and $v$ depend both on a \textit{flow} parameter $t$, and 
that $\bu(\ts{x},t)$ and $v(\ts{x},t)$ satisfy:
\begin{eqnarray} 
\Delta\bu_t&=&-\mbox{div}\left[\deru{h}(\det\nabla\bu-v)(\mbox{adj}
\nabla\bu)^T\right ], \label{GFpde1EF}\\
\Delta v_t&=&-\eps^2\Delta v-\deru{h}(\det\nabla\bu-v)
+\frac{1}{2\eps^2}\deru{\phi_\tau}(v),\label{GFpde2EF}
\end{eqnarray}
in $\mc{B}$ for $t>0$, and $\bu(\ts{x},t)=\lambda \ts{I}$ and $v(\ts{x},t)=0$ 
for all $\ts{x}\in\partial\mc{B}$ and $t\ge0$. The gradient flow equation leads 
to a descent method for
the minimization of \eqref{modfunctPF} over \eqref{set:reg}. After 
discretizing the partial derivatives with respect to ``$t$", one can use a
finite element method to solve the resulting flow equation. In particular, if 
we let $\Delta t>0$ be given, and set $t_{i+1}=t_i+\Delta t$ where $t_0=0$, we 
can approximate $\bu_t(\ts{x},t_i)$ and $v_t(\ts{x},t_i)$ with:
\[
\ts{z}_i(\ts{x})= \dfrac{\bu_{i+1}(\ts{x})-\bu_i(\ts{x})}{\Delta t},\quad
w_i(\ts{x})= \dfrac{v_{i+1}(\ts{x})-v_i(\ts{x})}{\Delta t},
\]
where $\bu_i(\ts{x})=\bu(\ts{x},t_i)$, etc. (We take $\bu_0(R)$ and $v_0$
to be some initial deformation and phase function satisfying the boundary 
conditions at $\partial\mc{B}$.) Inserting these approximations 
into the weak form of \eqref{GFpde1EF} and \eqref{GFpde2EF} together with the 
boundary conditions, and evaluating the right hand sides of the resulting 
equations at $\bu=\bu_i$ and $v=v_i$, we arrive at the following iterative
formula:
\begin{eqnarray}
&\displaystyle\int_\mc{B}\nabla\ts{z}_i\cdot\nabla\ts{q}\,\dif\ts{x}+ 
\int_\mc{B}\deru{h}(\det\nabla\bu_i-v_i)(\mbox{adj}\nabla\bu_i)^T\cdot\nabla\ts{
q} \,\dif\ts{x}=0,\label{floweq1}&\\
&\displaystyle\int_\mc{B}\nabla w_i\cdot\nabla r\,\dif\ts{x}+
\int_\mc{B}\left[\eps^2\nabla v_i\cdot\nabla r+ 
\left(-\deru{h}(\det\nabla\bu_i-v_i)
+\frac{1}{2\eps^2}\deru{\phi_\tau}(v_i)\right)r\right]\,\dif\ts{x},
\label{floweq2}&
\end{eqnarray}
for all functions $\ts{q}$ and $r$ sufficiently smooth such that 
$\ts{q}(\ts{x})=\ts{0}$ and $r(\ts{x})=0$ on $\partial\mc{B}$, and so that the 
integrals above are well defined. Given $\bu_i$ and $v_i$, one can solve the 
above equations for $\ts{z}_i$ and $w_i$ via some
finite element scheme, and then set
\begin{equation}\label{GFiter}
\ts{u}_{i+1}=\ts{u}_i+\Delta t\,\ts{z}_i,\quad 
v_{i+1}=v_i+\Delta t\,w_i,
\end{equation}
where the size of $\Delta t$ in each of these iterations is controlled 
individually. This process is repeated for $i=0,1,\ldots$, until both 
$\bu_{i+1}-\bu_i$ and $v_{i+1}-v_i$ are ``small'' enough, or some maximum value 
of ``$t$'' is reached, declaring the last computed solution pair as an 
approximate minimizer of \eqref{modfunctPF} over \eqref{set:reg}. In our 
simulations, equations \eqref{floweq1}--\eqref{floweq2} are solved in an 
alternate direction fashion, that is, given the pair $(\bu_i,v_i)$, we solve 
\eqref{floweq1} for $\bu_{i+1}$, then solve \eqref{floweq2} for $v_{i+1}$ using 
$(\bu_{i+1},v_i)$. These problems are solved using the finite element package 
freefem++ (see \cite{He2012}).

It follows from the results in \cite{Ho1992}, in the special case of an 
elastic fluid, that the minimizer of the 
corresponding two dimensional version of \eqref{rad_energy} over 
\eqref{rad_adm}, is given by
\begin{equation}\label{PF:radmin}
 r_c(R)=\sqrt{d_0R^2+\lambda^2-d_0},
\end{equation}
for $\lambda>\sqrt{d_0}$. Here $\deru{h}(d_0)=0$, where $h$ satisfies \eqref{h} 
and is strictly convex. In this case one can easily check that 
$\det\nabla\bu_c=d_0$, where $\bu_c=(r_c(R)/R)\,\ts{x}$ and $R=\norm{\ts{x}}$. 
For $\lambda\le\sqrt{d_0}$, the minimizer is 
$r_h(R)=\lambda R$. Now for the function $h$ in \eqref{hmonotD} with the data 
in \eqref{numdata}, we get that $d_0=1$. Thus $\lambda_c=1$ is the critical 
boundary displacement for radial cavitation. We want to emphasize, that in the 
simulations presented in this section, the numerical scheme works in a fully 
two dimensional setting, that is, not assuming radial symmetry for the 
functions. Moreover, the triangulations of the region $\mc{B}$ used in the 
finite element method, are not radially symmetric. For the purpose of 
presenting the results in this section, the computed finite element 
approximations were evaluated on a radially symmetric mesh. This data was then 
used in MATLAB\texttrademark~ to produce the figures. 

\begin{figure}
  \begin{tabular}{p{3.0in}p{4.0in}}
  \subfloat[Final approximation]
{\label{fig:mesh}\includegraphics[width=0.5\textwidth]{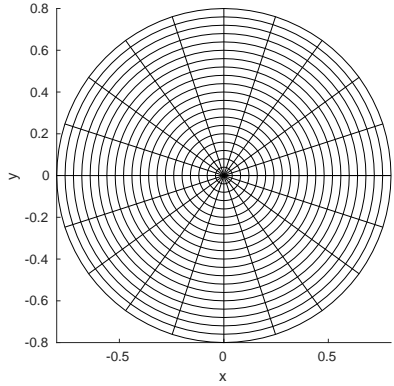}}&
  \subfloat[Determinant]
{\label{fig:fdef}{\includegraphics[width=0.5\textwidth]{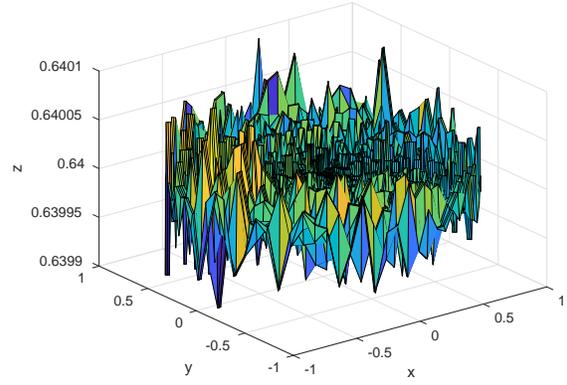}}}
\end{tabular}
\caption{Final approximation and its determinant computed by the decoupled 
method for the case $\lambda=0.8$.}\label{fig:1}
\end{figure}

\begin{figure}
\begin{center}
 \includegraphics[width=0.5\textwidth]{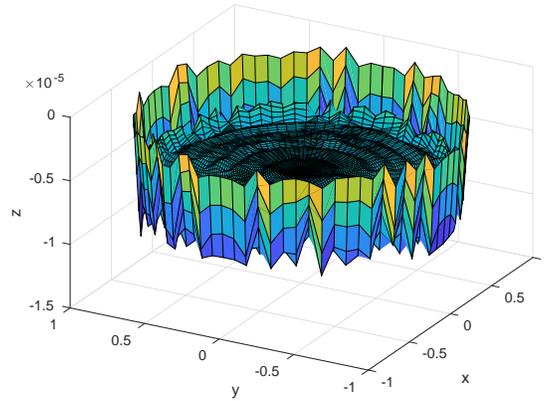}
 \end{center}
\caption{Final $v$ approximation computed by the decoupled 
method for the case $\lambda=0.8$.}\label{fig:2}
\end{figure}

In our first simulation, we consider the case of compression in which 
$\lambda=0.8$. Combining the results of the previous sections, it follows that 
the decoupled method converges to $\ts{u}_h(\ts{x})=\lambda\ts{x}$ for some 
sequences $\eps_j\To0$ and $\tau_j\To\infty$. To keep matters simple, we 
computed with fixed $\eps$ and $\tau$. In particular we used $\eps^2=10^{-2}$, 
$\tau=3$, a triangulation of $\mc{B}$ with 1750 triangles, and $M=10^4$ in 
\eqref{extphi}. We point out that $\bu_h$ and $v=0$ is not a pair solution of 
the Euler--Lagrange equations \eqref{pdeEF}. However as we will see now, with 
the selection of $\bu_0$ and $v_0$ as the finite element representations of 
$\bu_h$ and $v=0$ respectively, together with the selected parameters for the 
scheme mentioned above, the decoupled method produces an approximation pair 
very close to $\ts{u}_h$ and $v=0$. In Figure \ref{fig:1} we show the final 
computed deformation by the decoupled method together with its determinant. The 
final deformation is essentially a discretized version of $\ts{u}_h$, while the 
determinant ``oscillate'' around $0.64$ which is the determinant of gradient of 
$\bu_h$. In Figure \ref{fig:2} we show the computed $v$ which is of the order 
of $10^{-5}$, mostly negative, which helps to counter balance the two 
terms in \eqref{pde2EF} not including the laplacian. These 
approximations required slightly over 200 iterations of the 
gradient flow iteration, with an absolute error in \eqref{GFiter} of the order 
of $10^{-4}$. The computed energy by the decoupled method was $8.95075$ which 
compares well with the exact one which is $\pi h(0.64)=8.95670$, correct to the 
shown digits.

\begin{figure}
\begin{center}
 \includegraphics[width=0.7\textwidth]{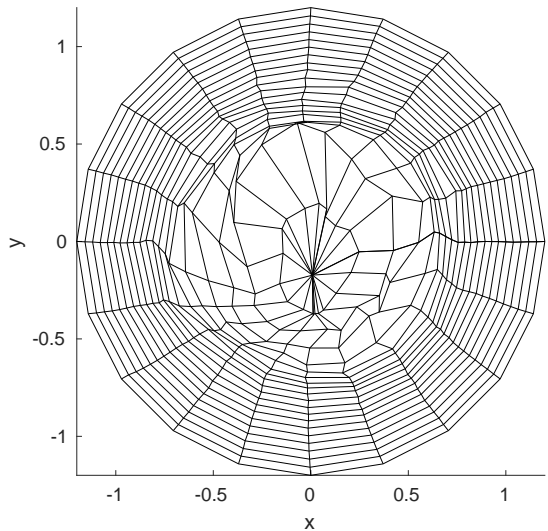}
 \end{center}
\caption{Computed approximation to $\bu_c$ by the decoupled 
method for the case $\lambda=1.2$.}\label{fig:3}
\end{figure}

\begin{figure}
  \begin{tabular}{p{3.0in}p{4.0in}}
  \subfloat[Norm of affine deformation]
{\label{fig:contafineexpan}\includegraphics[width=0.5\textwidth]{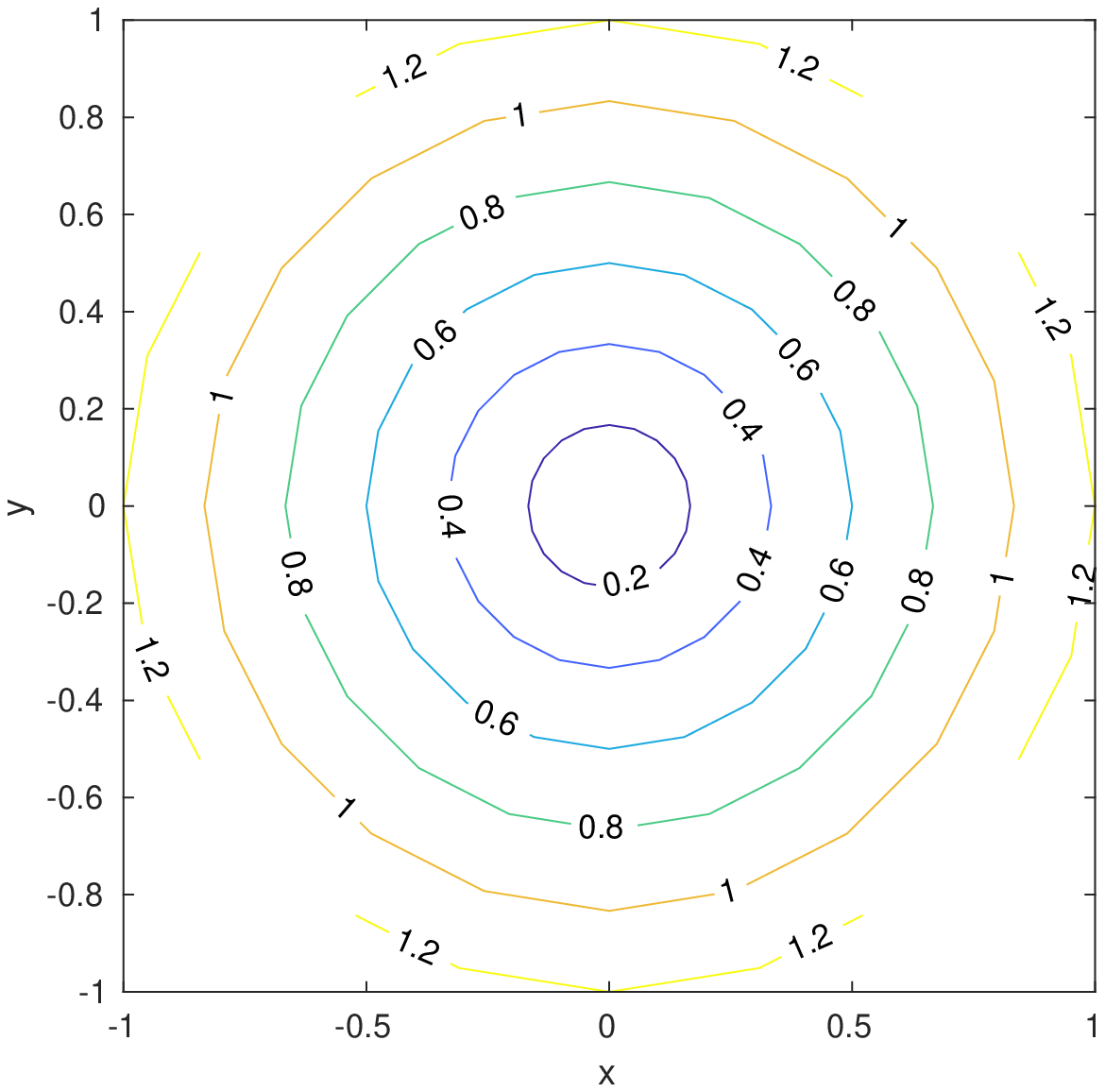}
} &
  \subfloat[Norm of final deformation]
{\label{fig:contdefexpan}{\includegraphics[width=0.5\textwidth]{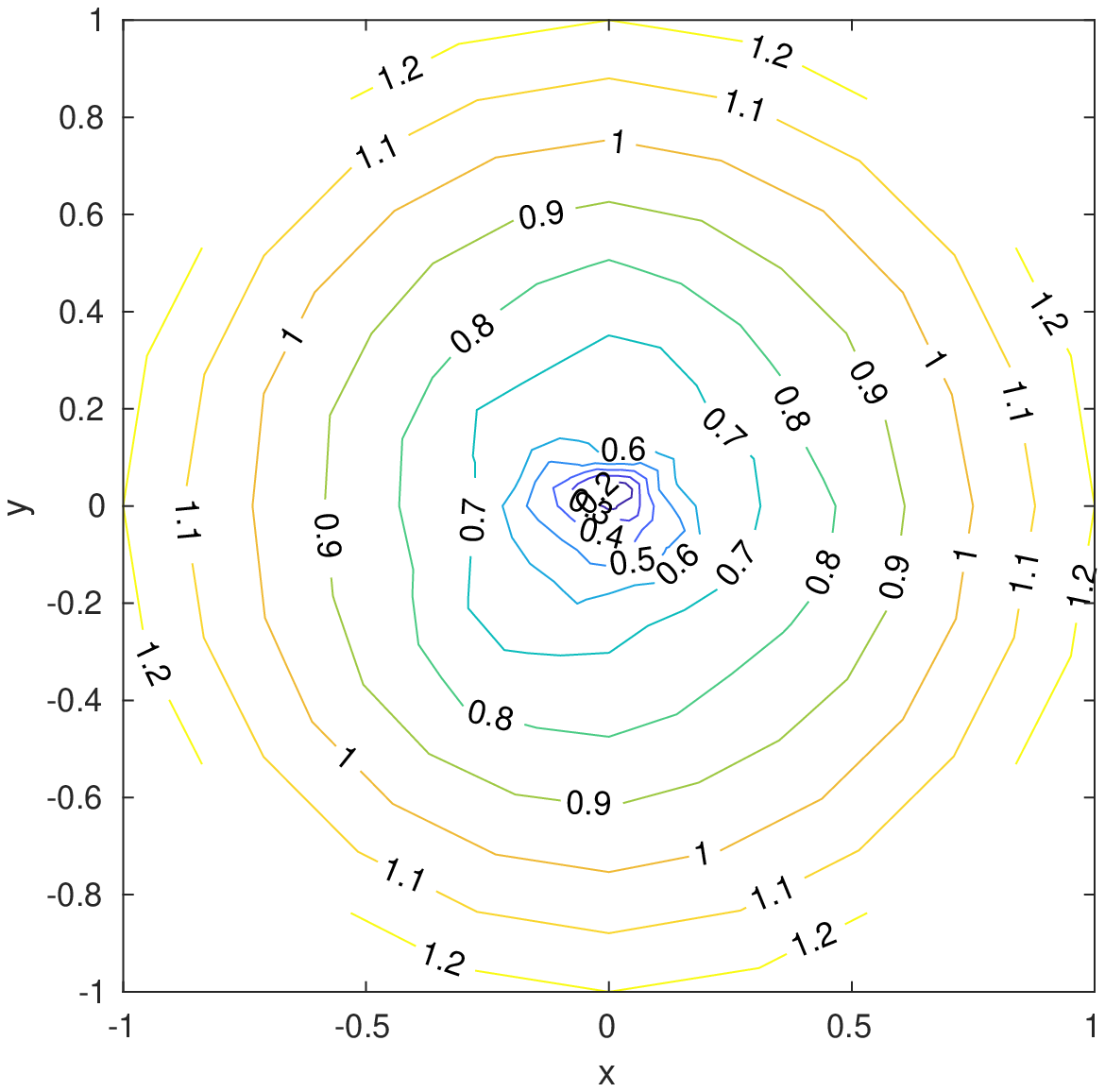}}}
\end{tabular}
\caption{Contour plots for $\norm{\lambda\ts{x}}$ and $\norm{\bu_\eps}$, for the 
case $\lambda=1.2$.}\label{fig:6}
\end{figure}

Our next simulation is for $\lambda=1.2$. In this case the decoupled method 
converges to $\bu_c=(r_c(R)/R)\,\ts{x}$ with $r_c(R)=\sqrt{R^2+0.44}$, and 
$\det\nabla\bu_c=1$. This deformation produces a cavity at the center 
$r_c(0)\approx0.6633$. For this simulation the parameters for the scheme were 
chosen as $\eps^2=10^{-3}$, $\tau=1$, a triangulation of $\mc{B}$ with 4534 
triangles, and $M=10^2$ in \eqref{extphi}. For $\bu_0$ we took the finite 
element version of $\bu_h$, and for $v_0$, the finite element version of 
$1.2\,\me^{-40(x^2+y^2)}$. We ran the gradient flow iteration for 30000 
iterations, with an absolute error in \eqref{GFiter} of the order 
of $10^{-5}$. The computed energy by the decoupled method was $6.29509$ 
which is comparable with the exact one which is $\pi h(1)=6.28318$, 
correct to the shown digits. In Figure \ref{fig:3} we show the computed 
computed approximation to $\bu_c$, which in what follows, we call $\bu_\eps$. 
The ``clearing'' in the central region in this picture, corresponds to material 
that has been expanded as to make the computed approximation $\bu_\eps$ to 
resemble the cavitating deformation $\bu_c$. From the figure we see that the 
approximate radius of this region is around $0.6$ which is consistent with 
cavity radius produced by $\bu_c$. Another way to measure the displacement of 
material from the central region in the reference configuration, is by looking 
at the contour graph of $\norm{\bu_\eps}$. In Figure \ref{fig:6} we show both, 
the contour plot of the norm of the affine deformation $\lambda\ts{x}$, and 
that for $\norm{\bu_\eps}$. This picture clearly shows how the contour curves 
from the reference configuration, get compressed very close to the origin, a 
consequence of the large displaments in this region. Note that the computed 
approximation is not radially symmetric, except for the region outside of the 
inner clearing. We think this lack of symmetry in the interior is due to the 
large deformations taking place to mimic the cavity by $\bu_c$, combined with 
the lack of symmetry of the finite element triangulation. 

\begin{figure}
  \begin{tabular}{p{3.0in}p{4.0in}}
  \subfloat[$\det\nabla\bu_\eps$]
{\label{fig:detexpan}\includegraphics[width=0.5\textwidth]{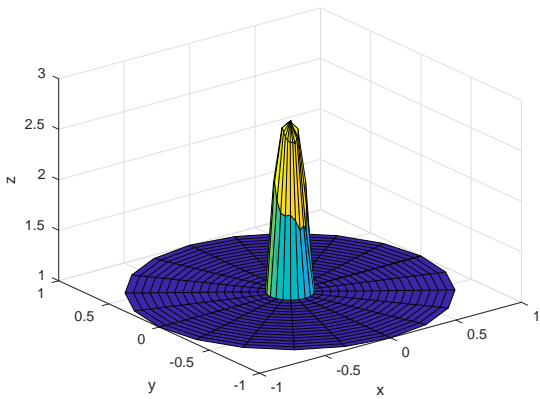}}&
  \subfloat[$\det\nabla\bu_\eps-v_\eps$]
{\label{fig:detmvexpan}{\includegraphics[width=0.5\textwidth]{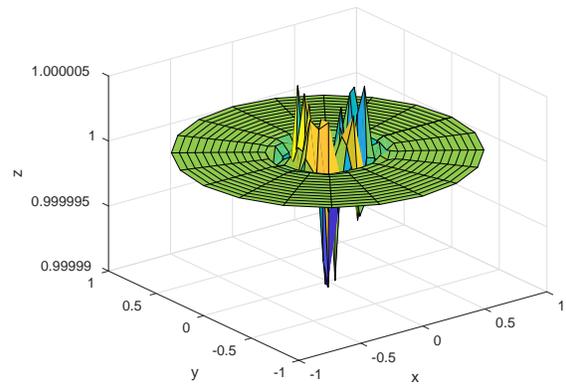}}}
\end{tabular}
\caption{Determinant of final approximation and the difference between this 
determinant with the computed $v$ approximation, for the case 
$\lambda=1.2$.}\label{fig:4}
\end{figure}

We show now in Figure \ref{fig:4} the 
graphs $\det\nabla\bu_\eps$ and $\det\nabla\bu_\eps-v_\eps$, where $v_\eps$ is 
the computed $v$ approximation. We note that the determinant of the gradient of 
$\bu_\eps$ is rather large close to the origin, indicative of the 
large deformations taking place in this region, while away from the center is 
essentially one. However the difference $\det\nabla\bu_\eps-v_\eps$ is very 
close to one globally, thus minimizing the $h$ term in the functional 
\eqref{modfunctPF}. As for the computed approximation $v_\eps$, we show in 
Figure \ref{fig:4} the initial $v$ approximation together with $v_\eps$. We can 
clearly see that the $v_\eps$ is concentrating around the origin, becoming 
larger in this region.

\begin{figure}
  \begin{tabular}{p{3.0in}p{4.0in}}
  \subfloat[Initial $v$]
{\label{fig:initvexpan}\includegraphics[width=0.5\textwidth]{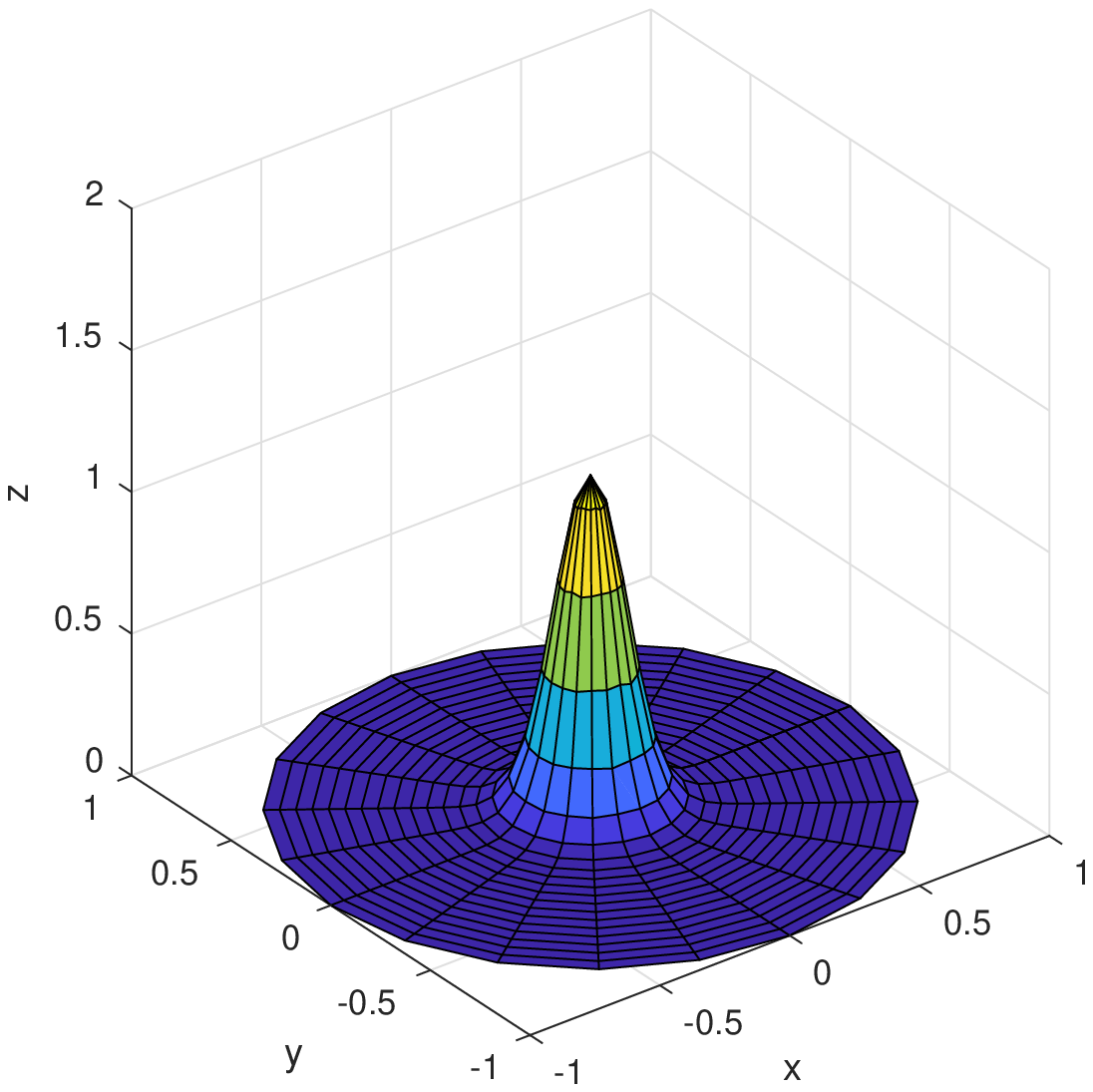}}&
  \subfloat[Final $v$]
{\label{fig:finalvexpan}
{\includegraphics[width=0.5\textwidth]{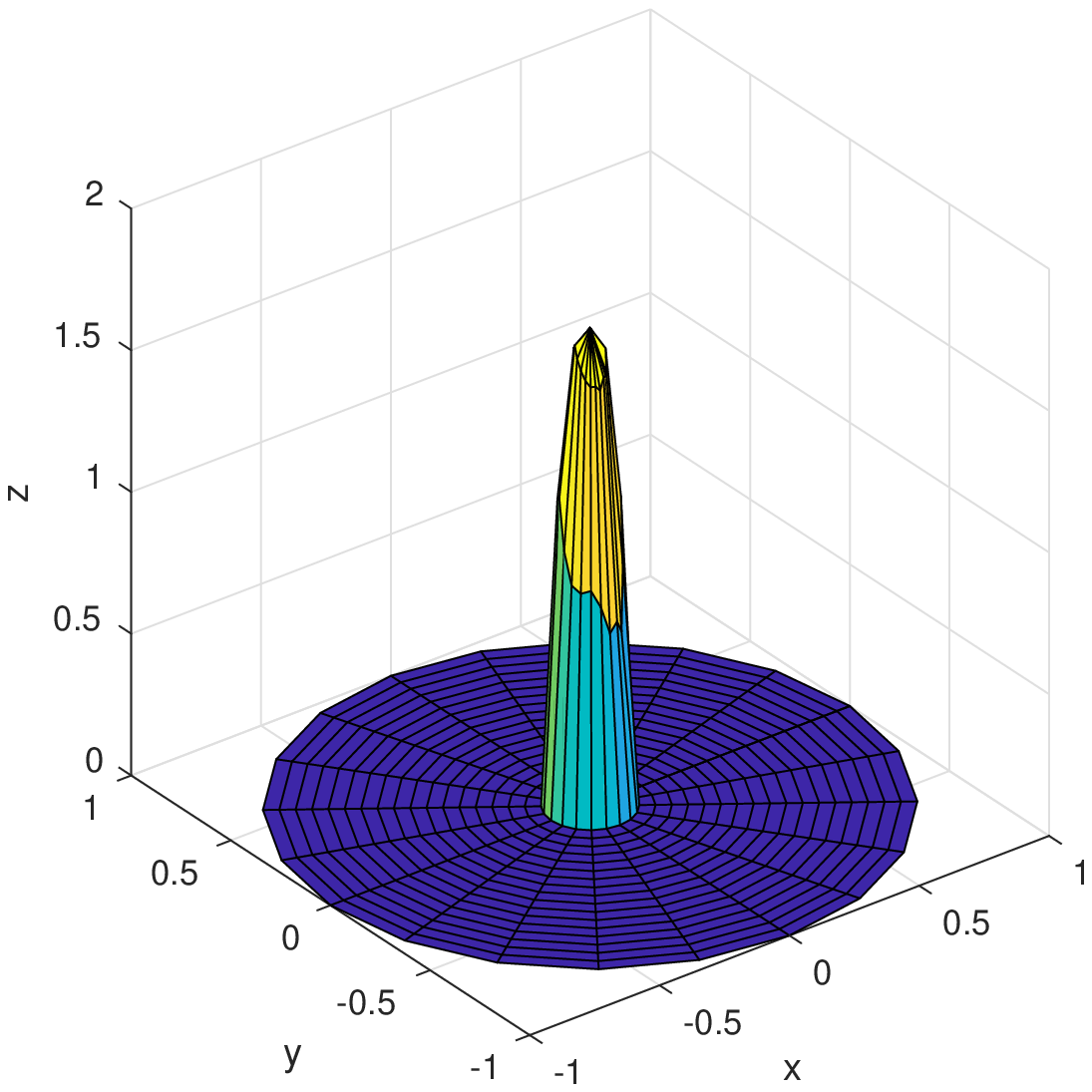}}}
\end{tabular}
\caption{Initial $v$ and the computed $v$ approximation, for the 
case $\lambda=1.2$.}\label{fig:5}
\end{figure}

\section{Concluding Remarks}
From the proof of Theorem \ref{thm:repprop} it becomes clear that the 
critical term in the stored energy function, in relation to the repulsion 
property, is the compressibility term, i.e., the function $h(\cdot)$ in 
\eqref{SEF}. This result is the main idea behind the method proposed in 
Section \ref{sec:num} and might explain why previous numerical schemes, such
as the element removal method developed by Li and coworkers (see, e.g.,
\cite{Li95})  or the use of ``punctured domains'' (see, e.g., 
\cite{SiSpTi2006}), have been successful.

As a practical matter, we mention that the numerical routine that one employs 
to solve the discrete versions of the minimization of \eqref{modfunct} over 
\eqref{set:reg}, must be ``aggressive'' enough, specially during the early 
stages of the minimization, to allow for actual increases in the 
intermediate approximate energies, which rules out the use of strictly 
descend methods. The reason for this is that, when needed, the scheme has to 
increase the 
phase function $v$ in regions where the determinant of the deformation gradient 
might become large. To do so, it might be necessary to overcome the maximum of 
the 
penalty function $\phi$ in \eqref{modfunct} resulting in an increase in the 
computed energy. One could try to avoid this by taking initial candidates 
for $v$ large, but this requires identifying regions where this is to be done, 
which in turn presumes knowledge of the location of the singularities. Although 
in general one can not assume such knowledge, it might be the case if the 
locations of possible flows in the material are known before hand.

Finally we did not address the question of the convergence of the minimizers of 
the discretized versions of \eqref{modfunct} over \eqref{set:reg}. Also we need 
to test the method on more general problems, other than the elastic fluid 
case. In particular, those problems in which the Lavrentiev phenomena 
takes place for boundary value problems in two dimensional elasticity among 
admissible continuous deformations. (See \cite{FoHrMi2003}.) These questions 
shall be pursued elsewhere.

\end{document}